\input amstex
\magnification=\magstep1
\input epsf
\input amssym.def
\input amssym
\pageno=1
\baselineskip 14 pt
\def \pop#1{\vskip#1 \baselineskip}

\font\gr=cmbx12

\def \et{\operatorname {et}}

\def \Fr{\operatorname {Fr}}
\def \Spec{\operatorname {Spec}}
\def \Spf{\operatorname {Spf}}
\def \lim{\operatorname {lim}}

\def \Sp  {\operatorname {Sp  }}

\def \fppf{\operatorname {fppf }}

\def \Spf{\operatorname {Spf}}
\def \sp{\operatorname {sp}}
\def \nsp{\operatorname {nsp}}

\def \Aut{\operatorname {Aut}}

\def \Vert{\operatorname {Vert}}
\def \Fr{\operatorname {Fr}}

\def \kum{\operatorname {kum}}

\def \Deg{\operatorname {Deg}}
\def \SS-Deg{\operatorname {SS-Deg}}
\def \SNS-Deg{\operatorname {SNS-Deg}}
 \def \DS-Deg{\operatorname {DS-Deg}}
\def \DNS-Deg{\operatorname {DNS-Deg}}

\pop {4}
\par
\noindent                                          
\centerline {\bf \gr Galois covers of degree $p$ : semi-stable}
\par
\noindent                                          
\centerline {\bf \gr reduction and Galois action}

\pop {3}
\noindent                                          
\centerline {\bf \gr Mohamed Sa\"\i di}

\pop {4}
\par
\noindent
{\bf \gr 0. Abstract.}\rm \ In this paper we study the semi-stable 
reduction of Galois covers of degree $p$ above semi-stable curves over a 
complete discrete valuation ring of inequal chracteristics $(0,p)$. We are 
also able to describe the Galois action on these covers in terms of some 
geometric and combinatorial datas in characteristic $p$ endowed with the 
action of the Galois group of the residue field.

\pop {2}
\par
\noindent
{\bf \gr 0. Introduction.}\rm \ In this paper we study the semi-stable 
reduction of Galois covers of degree $p$ above semi-stable curves over a 
complete discrete valuation ring of inequal chracteristics $(0,p)$, as well
as the Galois action on these datas. 
\pop {.5}
\par
More precisely, let $R$ be a complete discrete valuation ring of inequal
chracteristics $(0,p)$, with fraction field $K$, and residue field $k$.
In the first part of the paper we consider the following local situation. Let
$f:\Cal Y\to \Cal X$ be a Galois cover of degree $p$ between formal germs
of $R$-curves at the closed points $y$ and $x$ respectively, 
with $\Cal Y$ normal, and $\Cal X$ semi-stable i.e. $\Cal X$ is the formal 
germ of an $R$-curve at a closed point $x$ which is either a smooth point 
or an ordinary double point. It follows from the theorem of semi-stable 
reduction for curves that the germ $\Cal Y$ admits potentially a semi-stable 
reduction i.e. there exists, after eventually a finite extension of $R$,
a semi-stable model $\Tilde {\Cal Y}\to \Cal Y$ which is ``essencially'' 
unique. In particular, the Galois group of the cover $f$ acts on 
$\Tilde {\Cal Y}$ and the quotient $\Tilde {\Cal X}$ of $\Tilde {\Cal Y}$ 
by this action is a semi-stable blow-up of $\Cal X$. We have a canonical 
morphism $\Tilde f:
\Tilde {\Cal Y}\to \Tilde {\Cal X}$ which is Galois of degree $p$.
To the above cover $f:\Cal Y\to \Cal X$ we associate canonically some 
{\it degeneration datas} which determine completely
the special fibre $\Tilde {\Cal Y}_k:=\Tilde {\Cal Y}\times _Rk$ of
$\Tilde {\Cal Y}$. These
degeneration datas consist of the (canonically marked) tree associated to 
the special fibre $\Tilde {\Cal X}_k:=\Tilde {\Cal X}\times _Rk$ of
$\Tilde {\Cal X}$, plus some geometric datas which consists of 
a torsor $f_i:V_i\to U_i$ under a finite and flat group scheme of rank $p$
above each irreducible marked component $U_i$ of $\Tilde {\Cal X}_k$ such 
that the ``conductors'' of these torsors at the marked points satisfy 
certain ``compatibility'' conditions (these are what we called in [Sa] 
Kummerian mixed torsors) (cf. 1.2, 1.3 for a more precise definition). 
For an illustration of this 
situation we refeer to the example 1.2.1. Let's denote by $\Deg_p$
the set of isomorphism classes of such degeneration datas which are defined 
over an algebraic closure $\overline k$ of $k$ (cf. 1.2.2, 1.3.1 for a precise
definition). We show that the set 
$\Deg_p$ is a $G_k$-set, where $G_k$ is the Galois group of the separable 
closure of $k$ contained in $\overline k$ (cf. 1.2.3). Let $\overline K$ be 
an algebraic closure of $K$, and let $\Spec L$ be the geometric generic  
point of $\Cal X$. The \'etale cohomology group $H^1_{\et}(\Spec L,\mu_p)$
classifies (pointed)-geometric Galois covers of $\Cal X$ of degree $p$.
Moreover our result which exhibit the above degeneration datas can be
translated as the existence of a canonical specialisation map $\Sp:
H^1_{\et}(\Spec L,\mu_p)\to \Deg_p$. Note that both
$H^1_{\et}(\Spec L,\mu_p)$ and $\Deg_p$ are $G_K$-sets, where $G_K$
is the Galois group of $\overline K$ over $K$, and $G_K$ acts on
$\Deg_p$ via its canonical quotient $G_k$.
Our main result in this part of the paper is the following:

\pop {1}
\par
\noindent
{\bf \gr Theorem (1.2.7, 1.3.2).}\rm \ {\sl The above specialisation map 
$\Sp:H^1_{\et}(\Spec L,\mu_p)\to \Deg_p$ is surjective and $G_K$-equivariant.}

\pop {.5}
\par
In other words the association of degeneration datas to a cover
$f:\Cal Y\to \Cal X$ as above is compatible with the action of 
the Galois group on covers and on degeneration datas, 
Moreover one has a realisation result for such 
degeneration datas (cf. example 1.2.5 which illustrates the realistaion 
of degeneration datas). The proof of the above result relies heavily on the 
results in [Sa] and [Sa-1]. We also use formal patching techniques \`a 
la Harbater. As an application of the above result we construct in 2.1.8 
an example of a covers $f:\Cal Y\to \Cal X$ as above where $\Cal X$ is 
smooth at the closed point $x$, and where the special fibre 
$\Cal Y_k:=\Cal Y\times _Rk$ of $\Cal Y$ is singular and unibranche at the 
closed point $y$, and such that the configuration of 
the special fibre of the semi-stable model $\Tilde {\Cal Y}$ of $\Cal Y$ is 
not a tree-like. This indeed answers a question raised by Lorenzini
whether such a situation can occur in the inequal characteristic case.

\pop {.5}
\par
In the second part of the paper we study the global situation of a 
Galois cover $f:Y\to X$ of degree $p$ above a proper and semi-stable 
$R$-curve $X$ with $Y$ normal. As above, and using the theorem of 
semi-stable reduction for curves (cf. [De-Mu]), one obtains after eventually
a finite extension of $R$ a semi-stable model $\Tilde Y\to Y$ of $Y$ 
and a Galois cover $f:\Tilde Y\to \Tilde X$ of degree $p$, where 
$\Tilde X$ is a semi-stable blow-up of $X$. We associate then 
canonically to the cover $f:Y\to X$ some ``degeneration datas'' which 
determine completely the special fibre $\Tilde Y_k:=\Tilde Y\times _Rk$ of
the semi-stable model $\Tilde Y$ of $Y$ (cf. 2.2). This consists of the graph 
associated to the semi-stable $k$-curve $X_k$, and a given ``mixed torsor''
above $X_k$, plus given local degeneration datas at the critical points 
of this mixed torsor (cf. 2.2.1 for a more precise definition). 
For an illustration of this we refeer to the example 2.2.3.
Let's denote by $\Deg_p(X_k)$ the set of isomorphism 
classes of such degeneration datas which are defined over $\overline k$.
We show that this set is a $G_k$-set in a canonical way. Let 
$\overline \eta$ be the geometric generic  
point of $X$. The \'etale cohomology group $H^1_{\et}(\overline \eta,\mu_p)$
classifies (pointed)-geometric Galois covers of $X$ of degree $p$.
Moreover our result which exhibit the above degeneration datas can be
translated as the existence of a canonical specialisation map $\Sp:
H^1_{\et}(\overline \eta,\mu_p)\to \Deg_p(X_k)$. Note that both
$H^1_{\et}(\overline \eta,\mu_p)$ and $\Deg_p(X_k)$ are $G_K$-sets, where $G_K$
is the Galois group of $\overline K$ over $K$, and $G_K$ acts on
$\Deg_p$ via its canonical quotient $G_k$.
Our main result in the second part of the paper is the following: 

\pop {1}
\par
\noindent
{\bf \gr Theorem (2.3.2).}\rm \ {\sl The above specialisation map 
$\Sp:H^1_{\et}(\overline \eta,\mu_p)\to \Deg_p(X_k)$ is $G_K$-equivariant.}

\pop {.5}
\par
The above specialisation map $\Sp$ can not be surjective in this global 
case. However we prove the following result of realisation of degeneration 
datas:

\pop {1}
\par
\noindent
{\bf \gr Theorem (2.3.1).}\rm \ {\sl Let $X_k$ be a proper and semi-stable 
$k$-curve. Let $R$ be a complete discrete valuation ring of inequal 
characteristics with residue field
$k$, and fractions field $K$. Suppose given a 
$\overline k$-degeneration data $\deg (X_k)$ of rank $p$ 
associated to $X_k$ (in other words suppose given an element of $\Deg_p(X_k)$).
Then there exists, after eventually a finite extension of $R$, 
a proper and semi-stable $R$-curve $\Tilde X$ with smooth generic fibre 
and a special fibre $\Tilde X_k:=\Tilde X\times_Rk$
isomorphic to $X_k$, and a Galois cover
$\Tilde f:\Tilde Y\to \Tilde X$, such that the degeneration data 
associated to $f$ is isomorphic to the given data $\deg (X_k)$.} 

\pop {.5}
\par
We refeer to the example 2.3.4 which illustrate the realisation of global 
degeneration datas. Finally, I believe the ideas presented in this paper
provide the framework to construct Hurwitz spaces over $\Bbb Z$ whose fibre 
in characteristic zero classifies Galois covers of degree $p$ plus extra datas.
Also in the case of Galois covers of degree $p$ of the projective line over 
a $p$-adic field these results should lead to an algorithm which compute 
the semi-stable reduction of these covers at least in the case where the 

number of branched points is smaller than $p$.

\pop {1}
\par
\noindent
{\bf \gr I. Semi-stable reduction of Galois covers of degree $\bold {p}$
above formal germs of curves with Galois action.}
\rm
\pop {.5}      
\par
\noindent
{\bf \gr 1.0.}\rm \ In what follows we use the following notations: 
$R$ is a complete discrete valuation ring of inequal 
characteristics, with residue characteristic $p>0$, and
which contains a primitive $p$-th root of unity $\zeta$. 
We denote by $K$ the fraction field of $R$, 
$\pi $ a uniformising parameter, and $k$ the residue field of $R$. 
Let  $\overline K$ be a fixed 
algebraic closure of $K$, and let $G_K$  
be the Galois group of $\overline K$ of $K$. 
Let  $\overline R$ be the integral closure of $R$ in $\overline K$ which
is a valuation ring, and let $\overline k$ be the residue field of 
$\overline R$ which is an algebraic closure of $k$. 
Let $G_k$ be the Galois group of a separable closure of $k$ contained in 
$\overline k$. We have the following
canonical exact sequence, where $I_K$ denotes the inertia subgroup:        
$$0\to I_K\to G_K\to G_k\to 0$$
\pop {.5}
\par
In this section we will consider a formal germ $\Cal X$ of a 
semi-stable $R$-curve at a closed point, a Galois cover $f:\Cal Y\to \Cal X$ 
with group $G=\Bbb Z/p\Bbb Z$, and we will study the semi-stable reduction 
of $\Cal Y$, as well as the action of $G_K$ on these datas. 

\pop {.5}
\par
\noindent
{\bf \gr 1.1.}\ Let ${\Cal X}:=\Spec \hat {\Cal O}_{X,x}$ be the 
formal germ of an $R$-curve $X$ at a closed point $x$, and let
$f:\Cal Y\to \Cal X$ be a Galois cover with 
group $G\simeq \Bbb Z/p\Bbb Z$, such that  $\Cal Y$ is normal and 
local. It follows 
then easily 
from the
theorem of semi-stable reduction for curves (cf. [De-Mu], as well as the 
compactification process in [Sa-1] 2.3), that after 
eventually a finite extension $R'$ of $R$ with residue field $k'$,
and fractions field $K'$, the formal germ $\Cal Y$ has a semi-stable 
reduction. More precisely there exists a birational and proper morphism 
$\Tilde f :\Tilde {\Cal Y}\to \Cal Y'$, where  $\Cal Y'$ is the 
normalisation of $\Cal Y\times _R R'$, such that 
$\Tilde {\Cal Y}_{K'}\simeq \Cal Y'_{K'}$, and 
the following conditions hold:
\pop {.5}
\par
\noindent
(i) The special fibre $\Tilde {\Cal Y}_{k'}:=\Tilde {\Cal Y}\times _
{\Spec R'}\Spec k'$ of $\Tilde {\Cal Y}$ is reduced.
\par
\noindent
(ii) $\Tilde {\Cal Y}_k$ has only ordinary double points as singularities.
\pop {.5}
\par
Moreover there exists such a semi-stable model $\Tilde f :\Tilde {\Cal Y}
\to \Cal Y'$ which is {\bf minimal} for the above properties. In 
particular the 
action of $G$ on $\Cal Y'$ extends to an action
on $\Tilde {\Cal Y}$. Let $\Tilde {\Cal X}$ be the quotient of 
$\Tilde {\Cal Y}$ by $G$, which is a semi-stable model of $\Cal X$. One has
the following commutative diagram:

$$
\CD
\Tilde {\Cal Y}     @> \Tilde f >> \Cal Y '\\
   @Vg VV                  @Vf' VV   \\
\Tilde {\Cal X} @>\Tilde g>>  {\Cal X'}
\endCD
$$
\par
With the same notations as above, one can moreover choose the semi-stable 
models $\Tilde {\Cal Y}$ and $\Tilde {\Cal X}$ such that the
set of points $B_{K'}:=\{x_{i,K'}\}_{1\le i\le r}$, consisting of the branch 
locus in the morphism $\Cal Y'_{K'}\to \Cal X'_{K'}$, 
which we may assume to be {\bf rational}, specialise in {\bf smooth} (resp. 
{\bf smooth distincts}) points of $\Cal X'_{k'}$. We may also, after eventually
a finite extension of $K$, suppose that the double points of
$\Tilde {\Cal X}_k$ are rational. Moreover one can choose such 
$\Tilde {\Cal X}$ and $\Tilde {\Cal Y}$ which are minimal for these 
properties. We will denote by $\Tilde f ^{\nsp}:\Tilde {\Cal Y}^{\nsp}
\to \Cal Y'$ (resp. $\Tilde f ^{\sp}:\Tilde {\Cal Y}^{\sp}
\to \Cal Y'$) the minimal semi-stable model of $f$ such that the points
of $B_K$ specialise in smooth (resp. smooth distincts) points of 
$\Tilde {\Cal X}^{\nsp}:=\Tilde {\Cal Y}^{\nsp}/G$ (resp. of 
$\Tilde {\Cal X}^{\sp}:=\Tilde {\Cal Y}^{\sp}/G$). We call 
$\Tilde {\Cal Y}^{\nsp}$ (resp. $\Tilde {\Cal Y}^{\sp}$) the minimal {\bf
non split} (resp. {\bf split}) semi-stable model of $\Cal Y$. We have the 
following commutative diagrams: 

$$
\CD
\Tilde {\Cal Y}^{\nsp}     @> \Tilde f ^{\nsp}>> \Cal Y' \\
   @Vg^{\nsp} VV                  @Vf' VV   \\
\Tilde {\Cal X}^{\nsp} @>\Tilde g^{\nsp}>>  {\Cal X'}
\endCD
$$

and

$$
\CD
\Tilde {\Cal Y}^{\sp}     @> \Tilde f^{\sp} >> \Cal Y' \\
   @Vg^{\sp} VV                  @Vf' VV   \\
\Tilde {\Cal X}^{\sp} @>\Tilde g^{\sp}>>  {\Cal X'}
\endCD
$$

\pop {.5}
\par    
Moreover the morphism $\Tilde f^{\sp}:\Tilde {\Cal Y}^{\sp}\to \Cal Y'$ 
factors throught $\Tilde f^{\nsp}$ as follows: 
$\Tilde {\Cal Y}^{\sp}\to \Tilde 
{\Cal Y}^{\nsp}\to \Cal Y'$, and the later commutative diagram factors 
through the first one. In the case where $\Cal X$ is the formal germ of a 
semi-stable $R$-curve at a closed point $x$, the 
fibre $(\Tilde {g}^{\nsp})^{-1}(x)$ (resp. 
$(\Tilde {g}^{\sp})^{-1}(x)$) of the closed point $x$ 
in $\Tilde {\Cal X}^{\nsp}$
(resp. in $\Tilde {\Cal X}^{\sp}$)
is a tree $\Gamma ^{\nsp}$ (resp. $\Gamma ^{\sp}$) of projective lines. 
This tree is canonically 
endowed with some ``degeneration datas'' that we will exhibit below, 
and which 
follow mainly from the results in [Sa] and [Sa-1].

\pop {.5}
\par
\noindent
{\bf \gr 1.2.}\rm \ We will use the same notations as in 1.1. We first 
consider the case where $\Cal X\simeq \Spf
A$ is the formal germ of a semi-stable $R$-curve at 
a {\bf smooth} point $x$. Let $R'$ be a finite extension of $R$ as in 1.1,
and let $\pi '$ be a uniformiser of $R'$. We assume that $\Cal X$ is 
geometrically connected, in which case $\Cal Y$ is 
also geometrically connected. Below we exhibit the {\bf degeneration datas}
associated to the non split (resp. split) semi-stable reduction of 
$\Cal Y$, and which are consequences of the results in [Sa] and [Sa-1].
\pop {.5}
\par
\noindent
{\bf \gr Deg.1.}\rm \ Let $\wp:=(\pi')$ be the ideal of 
$A':=A\otimes _R R'$ generated by $\pi'$, 
and let $\hat A'_{\wp}$ be the 
completion of the localisation of $A'$ at $\wp$. Let 
$\Cal X'_{\eta}:=\Spf \hat A'_{\wp}$ be the boundary of $\Cal X'$, and let 
$\Cal X'_{\eta}\to \Cal X'$ be the canonical 
morphism. Consider the following cartesian diagram:

$$
\CD
{\Cal  Y}_{\eta} @>f_{\eta}>>  {\Cal X}'_{\eta}  \\                         
   @VVV                @VVV   \\
Y'    @> f'>> \Cal X' \\
\endCD
$$

Then $f_{\eta}:\Cal Y_{\eta}\to \Cal X'_{\eta}$ is a torsor under a 
commutative finite and flat 
$R'$-group scheme $G_{R'}$ of rank $p$. Let $\delta$ be 
the degree of the different associated to the torsor $f_{\eta}$. One has the 
degeneration  type $(G_{k'},m,h)$ of the torsor $f_{\eta}$ (cf. [Sa] 3.2)
which is canonically associated to $f$. The geometric genus $g_y$ of the
point $y$ equals $(r-m-1)(p-1)/2$ (cf. [Sa-1] 3.1.1), where $r$ is the 
cardinality of $B_{K'}$.
\pop {.5}
\par
\noindent
{\bf \gr Deg.2.}\rm \ The fibre 
$(\Tilde {g}^{\nsp})^{-1}(x)$ (resp. 
$(\Tilde {g}^{\sp})^{-1}(x)$) of the closed point $x$ of $\Cal X'$ 
in $\Tilde {\Cal X}^{\nsp}$
(resp. in $\Tilde {\Cal X}^{\sp}$)
is a tree $\Gamma ^{\nsp}$ (resp. $\Gamma ^{\sp}$) of projective 
lines, $\Gamma ^{\nsp}$ is a sub-tree of $\Gamma ^{\sp}$. Let 
$\Vert (\Gamma ^{\nsp}):=\{X_i\}_{i\in I^{nsp}}$ (resp. 
$\Vert (\Gamma ^{\sp}):=\{X_i\}_{i\in I^{sp}}$) be the set of 
irreducible components of $(\Tilde {g}^{\nsp})^{-1}(x)$ (resp.
of $(\Tilde {g}^{\sp})^{-1}(x)$), which are also the 
vertices of the tree $\Gamma ^{\nsp}$ (resp. $\Gamma ^{\sp}$). 
The tree $\Gamma^{\sp}$ (hence also $\Gamma^{\nsp}$)
is canonically endowed with an origin vertex $X_{i_0}$,
which is the unique irreducible component of 
$(\Tilde {g}^{\sp})^{-1}(x)$
which meets the point $x$. We fix
an orientation of the tree $\Gamma^{\sp}$ starting from $X_{i_0}$ 
in the direction of the ends. Such an orientation induces of 
course an orientation of the subtree $\Gamma^{\nsp}$. 
\pop {.5}
\par
\noindent
{\bf \gr Deg.3.}\rm \ For each $i\in I^{\nsp}$ (resp. $I^{\sp}$), 
let $\{x_{i,j}\}_{j\in S_i}$ be the 
set of points of $X_i$ in which 
specialise some points of $B_{K'}$ ($S_i$ may be empty), say in 
each point $x_{i,j}$ specialise $r_{i,j}$
points of $B_K$. If $i\in I^{\sp}$, 
and $S_i$ is non empty, we have $r_{i,j}=1$. Also let  
$\{z_{i,j}\}_{j\in D_i}$ be the set of points of the 
irreducible component $X_i$ where $X_i$ meets the rest of the 
components of $\Tilde {\Cal X}^{\nsp}_{k'}$ (resp. of 
$\Tilde {\Cal X}^{\sp}_{k'}$). These are the double points of 
$\Tilde {\Cal X}^{\nsp}_{k'}$ (resp. $\Tilde {\Cal X}^{\sp}_{k'}$) 
supported by $X_i$. We denote 
by $B_{k'}^{\nsp}$ (resp. $B_{k'}^{\sp}$) the set of all points 
$\cup _{i\in I^{\nsp}}\{x_{i,j}\}_{j\in S_i}$
(resp. $\cup _{i\in I^{\sp}}\{x_{i,j}\}_{j\in S_i}$), which is the set of 
specialisation of the branch locus $B_{K'}$, and by $D_{k'}^{\nsp}$ 
(resp. $D_{k'}^{\sp}$) 
the set of double points of $\Tilde {\Cal X}^{\nsp}_{k'}$ (resp.
$\Tilde {\Cal X}^{\sp}_{k'}$).  In particular
$x_{i_0,j_0}:=x$ is a double point of $\Tilde {\Cal X}^{\nsp}_{k'}$ (resp.
of $\Tilde {\Cal X}^{\sp}_{k'}$). 
\pop {.5}
\par
\noindent
{\bf \gr Deg.4.}\rm \ Let $\Cal U^{\nsp}:=\Cal X^{\nsp}-\{B_{k'}^{\nsp}
\cup D_{k'}^{\nsp}\}$ (resp. $\Cal U^{\sp}:=\Cal X^{\sp}-\{B_{k'}^{\sp}
\cup D_{k'}^{\sp}\}$). The cover $f$ induces a 
{\it mixed torsor} $f^{\nsp}:\Cal V\to \Cal U^{\nsp}$ above 
$\Cal U^{\nsp}$ (resp. $f^{\sp}:\Cal V'\to \Cal U^{\sp}$ above 
$\Cal U^{\sp}$), whose special fibre
$f^{\nsp}_{k'}:\Cal V_{k'}\to \Cal U^{\nsp}_{k'}$ 
(resp. $f^{\sp}_{k'}:\Cal V'_{k'}\to \Cal U^{\sp}_{k'}$) is an element of 
$H^1_{\fppf}(\Cal U^{\nsp}_{k'})_p$
(resp. an element of $H^1_{\fppf}(\Cal U^{\sp}_{k'})_p$) (cf.
[Sa] 1.6 for the definition of $H^1_{\fppf}(\ )_p$). In particular 
the restriction of $f^{\nsp}$ (resp. $f^{\sp}$) to each connected component 
$\Cal U_i$ of $\Cal U^{\nsp}$ (resp. $\Cal U^{\sp}$) is a torsor 
$f_i:\Cal V_i\to \Cal U_i$ under a commutative finite and flat $R'$-group 
scheme $G_{R',i}$ of rank $p$, and 
$f_{i,k'}:\Cal V_{i,k'}\to \Cal U_{i,k'}:={\Cal U}_i\times _{R'}{k'}$ 
is a torsor under the $k'$-group scheme $G_{k',i}:=G_{R',i}\times_{R'} 
{k'}$. Moreover if $G_{k',i}$ is radicial, and if $\omega _i$ is the 
associated differential form, then 
the set of zeros and poles of $w_i$ is necessarily contained in  
$\{x_{i,j}\}_{j\in S_i}\cup \{z_{i,j}\}_{j\in  D_i}$, as 
$\Cal V_{i,k'}$ is smooth (cf. [Sa], I).
\pop {.5}
\par
\noindent
{\bf \gr Deg.5.}\rm \ Each smooth point $x_{i,j}\in B_{k'}^{\nsp}$ 
(resp. $\in B_{k'}^{\sp}$) is endowed via $f$ with degeneration datas 
on the boundary of the formal fibre at $x_{i,j}$, as in Deg.1 above,  
and which satisfy certain {\it compatibility  conditions}. 
More precisely for each smooth point 
$x_{i,j}$ we have the reduction type $(G_{k',i},m_{i,j},h_{i,j})$
on the boundary of the formal fibre at this point, induced by $g^{\nsp}$
(resp. $g^{\sp}$), and such that $r_{i,j}=m_{i,j}+1$, as a consequence of 
[Sa-1] 3.1.1. Also if $x_{i,j}\in B_{k'}^{\sp}$ (in other words if 
$i\in I^{\sp}$, and $S_i$ is non empty), then necessarily 
$G_{k',i}=\mu _p$, $m_{i,j}=-1$ and $h_{i,j}=0$ (cf. [Sa-1, 3.1.2).
\pop {.5}
\par  
\noindent
{\bf \gr Deg.6.}\rm \ Each double point
$z_{i,j}=z_{i',j'}\in X_i\cap X_{i'}$ of $\Tilde X^{\nsp}$ (resp. 
$\Tilde X^{\sp}$) is endowed with degeneration datas 
$(G_{i,k'},m_{i,j},h_{i,j})$ and $(G_{i',k'},m_{i',j'},h_{i',j'})$, induced 
by $g^{\nsp}$ (resp. $g^{\sp}$), on the two boundaries of the formal 
fibre at this point as in Deg.1, such that $m_{i,j}+m_{i',j'}=0$, and 
$h_{i,j}+h_{i',j'}=0$, as a consequence of [Sa-1] 3.2.3. In particular
$m_{i_0,j_0}+m=0$, and $h_{i_0,j_0}+h=0$. In other words
the above induced by $f$ element of $H^1_{\fppf}(\Cal U^{\nsp}_{k'})_p$
(resp. of $H^1_{\fppf}(\Cal U^{\sp}_{k'})_p$) is indeed an element of
$H^1_{\fppf}(\Cal U^{\nsp}_{k'})^{\kum}_p$
(resp. an element of $H^1_{\fppf}(\Cal U^{\sp}_{k'})^{\kum}_p$) (cf.
[Sa] 1.7 for the definition of $H^1_{\fppf}(\ )^{\kum})$.
\pop {.5}
\par
\noindent
{\bf \gr Deg.7.}\rm \ For each double point
$z_{i,j}=z_{i',j'}\in X_i\cap X_{i'}$ of $\Tilde X^{\nsp}$ (resp. 
$\Tilde X^{\sp}$), let $e_{i,j}$ be the thikeness of $z_{i,j}$.
Then $e_{i,j}=pt_{i,j}$ is necessarily divisible by $p$. For each 
irreducible component
$X_i$, ${i\in I^{nsp}}$ (resp. $i\in I^{sp}$), let $\eta _i$
be the generic point of $X_i$. Let $\delta _i$ be the degree of different
above $\eta _i$ in the cover $g^{\nsp}$ (resp.  $g^{\sp}$). Then 
we have $\delta _i-\delta _{i'}=t_{i,j}m_{i,j}(p-1)$ as follows from 
[Sa-1] 3.2.5.
In particular $\delta -\delta _{i_0}=t_{i_0,j_0}m(p-1)$,
and $\vert D_i\vert \delta _i=\sum _{j\in D_i}(\delta _j+t_{i,j}m_{i,j}(p-1))$,
where $\vert D_i\vert$ denotes the cardinality of $D_i$.
\pop {.5}
\par
\noindent
{\bf \gr Deg.8.}\rm \ The contribution to the 
arithmetic genus $g_y$ of the point $y$ was 
computed in [Sa-1] 3.1.1, in terms of the 
degeneration data $(G_{k'},m,h)$ on the boundary of the formal fibre at $x$,
and the cardinality $r$ of $B_K$. It follows also from the above 
considerations and after easy calculation that: 
$g_y=\sum _{i\in I_{\et}}(-2+\sum _{j\in S_i}(m_{i,j}+1)
+\sum _{j\in D_i}(m_{i,j}+1))(p-1)/2$,   
where $I_{\et}$ denote the subset of $I^{\nsp}$, or $I^{\sp}$, which
among to the same, consisting of those $i$ for which
$G_{k',i}$ is \'etale.
\pop {.5}
\par
\noindent
{\bf \gr 1.2.1. Example.}\rm \ In the following we give 
an example where one can 
exhibit the degeneration datas associated to a Galois cover $f:\Cal Y\to
\Cal X$ of degree $p$ where $\Cal X\simeq \Spf R[[T]]$ is the formal germ 
of a smooth point. More precisely, for $m>0$ an integer prime to $p$ 
consider the cover given generically by the equation 
$X^p=1+\lambda ^p(T^{-m}+\pi T^{-m+1})$ (this is the 
example 1 in [Sa-1], 3.2.4, with $m'=m+1$). Hier $r=m+2$ and this cover has a 
reduction of type $(\Bbb Z/p\Bbb Z,m,0)$ on the boundary. In particular the 
geometric genus $g_y$ of the closed point $y$ of $\Cal Y$ 
equals $(p-1)/2$. The non 
split degeneration data associated to the above cover consists necessarily of
a treee with only one vertex and no edges, i.e. a unique projective line
$X_1$ with a marked $\overline k$-point $x_1$, and an \'etale torsor 
$f_1:V_1\to U_1:=X_1-\{x_1\}$ above $U_1$ with conductor $2$ at 
the point $x_1$. 
\pop {1}
\epsfysize=3.4cm
\centerline{\epsfbox{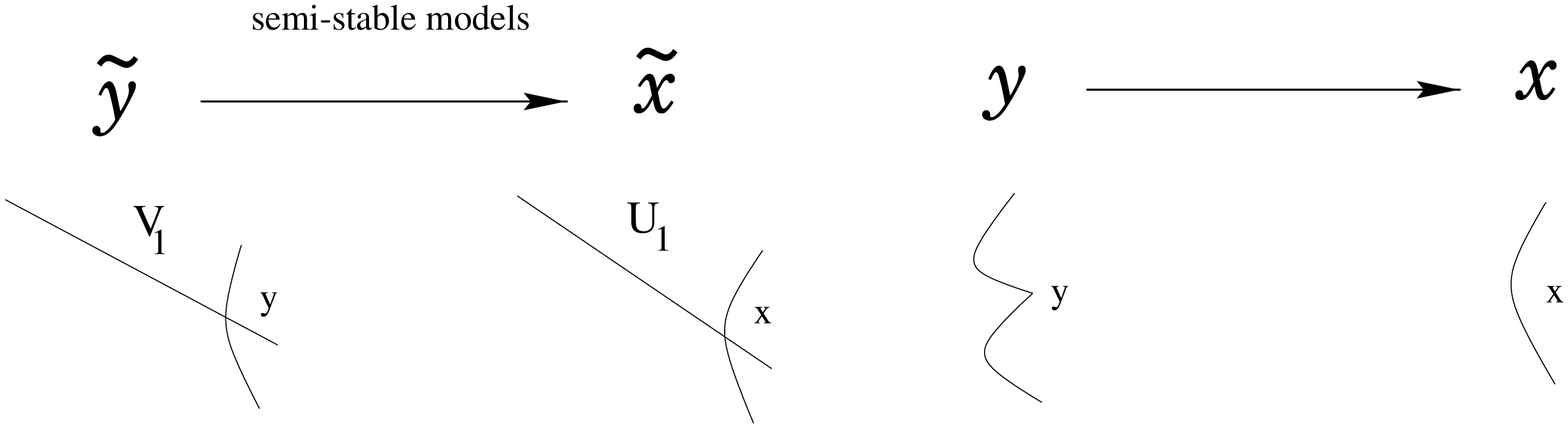}}

\par
The above considerations leads naturally to the following abstract
geometric and combinatorial definition of degeneration datas.
\pop {.5}
\par
\noindent
{\bf \gr 1.2.2. Definition.}\rm \ Let $k'$ be an algebraic extension 
of $k$. A $k'$-{\bf simple non split (resp. split) 
degeneration data} Deg(x) of type $(r,(G_{k'},m,h))$ and rank $p$
consists of the following datas:
\pop {.5}
\par
\noindent
{\bf \gr Deg.1.}\rm \ $G_{k'}$ is a commutative finite and 
flat $k'$-group scheme of rank $p$, $r\ge 0$ is an integer, $m$ is an 
integer prime to $p$ such that $r-m-1\ge 0$, and $h\in \Bbb F_p$ equals 
$0$ unless $G_{k'}=\mu _p$ and $m=0$.

\pop {.5}
\par
\noindent
{\bf \gr Deg.2.}\rm\  $\Gamma:=X_{k'}$ is an oriented tree of projective 
lines with vertices $\Vert (\Gamma):=\{X_i\}_{i\in I}$ (resp. 
$\Vert (\Gamma):=\{X_i\}_{i\in I^{\sp}}$), endowed 
with an origin vertex $X_{i_0}$, and a marked smooth point $x:=x_{i_0,j_0}$ 
on $X_{i_0}$ which is a $k'$-rational point. We denote by 
$\{z_{i,j}\}_{j\in D_i}$ the set of double points, or (non oriented) 
edges, of $\Gamma$ which are supported by $X_i$, and which we assume 
to be $k'$-rational points.

\pop {.5}
\par
\noindent
{\bf \gr Deg.3.}\rm \ For each vertex $X_i$ of $\Gamma$ is given a set, 
may be empty, of smooth marked $k'$-rational points $\{x_{i,j}\}_{j\in S_i}$.

\pop {.5}
\par
We will summarise the datas Deg.2 and Deg.3 by saying that we are given 
a $k'$-{\bf rational marked tree} of projective lines, or a
marked semi-stable $k'$-curve with arithmetic genus $0$.

\pop {.5}
\par
\noindent
{\bf \gr Deg.4.}\rm \ For each $i\in I$ (resp. $i\in I^{\sp}$), is given a 
torsor $f_{i,k'}:V_{i,k'}\to U_{i,k'}:= 
X_{i,k'}-\{\{x_{i,j}\}_{j\in S_i}
\cup \{z_{i,j}\}_{j\in D_i}\}$ under a commutative finite and flat 
$k'$-group scheme $G_{k',i}$ of rank $p$, where
$X_{i,k'}:=X_i\times _{\Spec k}\Spec {k'}$, 
with $V_{i,k'}$ {\bf smooth}. Let $U_{k'}$ be the open subscheme
of $X_{k'}$ obtained by deleting the double and marked points
of $X_{k'}$. Then the above data determines an element of
$H^1_{\fppf}(U_{k'})_p$ (cf. [Sa], I). Moreover if $G_{k',i}$ is 
radicial, and if $\omega _i$ is the associated differential form, then 
the set of zeros and poles of $w_i$ is necessarily contained in  
$\{x_{i,j}\}_{j\in S_i}\cup \{z_{i,j}\}_{j\in D_i}$, as
$V_{i,k'}$ is smooth. Also if $i\in I^{\sp}$, and $S_i$ is non empty, we 
assume that $G_{i,k'}=\mu_p$.

\pop {.5}
\par
\noindent
{\bf \gr Deg.5.}\rm \ For each $i\in I$ (resp. $i\in I^{\sp}$), are given for 
each $j\in S_i$ the integres $(m_{i,j},h_{i,j})$, where $m_{i,j}$
is the conductor of the torsor $f_i$ at the point $x_{i,j}$, 
and $h_{i,j}$ its residue at this point. Let $r_{i,j}:=m_{i,j}+1$.
If $i\in I^{\sp}$ then we assume that $r_{i,j}=1$. Also for the marked 
point $x_{i_0,j_0}$ is given the data $(m_{i_0,j_0},h_{i_0,j_0})$, such that
$m_{i_0,j_0}+m=0$, and $h_{i_0,j_0}+h=0$.

\pop {.5}
\par
\noindent
{\bf \gr Deg.6.}\rm \ For each double point 
$z_{i,j}=z_{i',j'}\in X_i\cap X_{i'}$ is given the pair $(m_{i,j},h_{i,j})$
(resp. $(m_{i',j'},h_{i',j'})$), where $m_{i,j}$ (resp. $m_{i',j'}$)
is the conductor of the torsor $f_{i,k'}$ (resp. $f_{i',k'}$) at the point 
$z_{i,j}$ (resp. $z_{i',j'}$), and $h_{i,j}$ (resp. $h_{i',j'}$) its 
residue at this point. These datas must satisfy $m_{i,j}+m_{i',j'}=0$, 
and $h_{i,j}+h_{i',j'}=0$. In other words the element
of $H^1_{\fppf}(U_{k'})_p$ determined by the data Deg.4 must belong
to the subgroup $H^1_{\fppf}(U_k)_p^{\kum}$ (cf. [Sa]. I).

\pop {.5}
\par
\noindent
{\bf \gr Deg.7.}\rm \ For each irreducible component
$X_i$ of $\Gamma$ is given an integer $\delta _i\le v_K(p)$ which
is divisible by $p-1$. For each double point $z_{i,j}=z_{i',j'}\in 
X_i\cap X_{i'}$ of $\Gamma$ is given an integer 
$e_{i,j}=pt_{i,j}$ divisible by $p$, such that with the same notations as 
above we have $\delta _i-\delta _{i'}=m_{i,j}t_{i,j}(p-1)$.

\pop {.5}
\par
\noindent
{\bf \gr Deg.8.}\rm \  Let $I_{\et}$ denote the subset of $I$ 
(resp. of $I^{\sp}$) consisting of those $i$ for which
$G_{k',i}$ is \'etale. Then the following equality should hold: 
$(r-m-1)(p-1)/2
=\sum _{i\in I_{\et}}(-2+\sum _{j\in S_i}
(m_{i,j}+1)+\sum _{j\in D_i}(m_{i,j}+1))(p-1)/2$. The integer 
$g_x:=(r-m-1)(p-1)/2$ is called the {\bf genus} of the degeneration 
data Deg(x).

\pop {.5}
\par
There is a natural notion of isomorphism of simple split (resp. non split) 
degeneration datas of a given type over a given algebraic extension of $k$. 
We will denote by $\bold {\SS-Deg_p}$ (resp. $\bold {\SNS-Deg_p}$) the set of 
isomorphism classes of $\overline k$-simple split (resp. 
$k'$-simple non split) degeneration data of rank $p$.

\pop {.5}
\par
\noindent
{\bf \gr 1.2.3. Galois action on degeneration datas.}\rm 
\pop {.5}
\par 
The Galois group
$G_k$ of the separable closure of $k$ contained in $\overline k$ 
acts in a canonical way on the
set ${\SS-Deg_p}$ (resp. ${\SNS-Deg_p}$)
of isomorphism classes of $\overline k$-simple split (resp. non split)
degeneration datas of rank $p$, in a way that is compatible with the 
group law on $G_k$. In other words  ${\SS-Deg_p}$ 
and $\SNS-Deg_p$ are $G_k$-sets. 
More precisely, to a $\overline k$-simple split 
(resp. non split) degeneration data Deg(x) of type 
$(r,(G_{\overline k},m,h))$, and an element $\sigma\in G_k$, 
we associate the simple degeneration data $\Deg(x)^{\sigma}$ of the same type 
$(r,G_{\overline k},m,h))$. We explain briefly the defining datas for
$\Deg(x)^{\sigma}$:

\pop {.5}
\par
\noindent
{\bf \gr Deg$^\bold {\sigma}$.1.}\rm \ The datas Deg$^{\sigma}$.1
 are the same as in Deg.1.

\pop {.5}
\par
\noindent
{\bf \gr Deg$^\bold {\sigma}$.2. and Deg$^{\sigma}$.3.}\rm\ 
The $\overline k$-rational marked tree $\Gamma^{\sigma}$ is the marked tree 
$\Gamma$.

\pop {.5}
\par
\noindent
{\bf \gr Deg$^\bold {\sigma}$.4.}\rm \ For each $i\in I$ 
(resp. $i\in I^{\sp}$), 
is given the torsor $f_i^{\sigma}:\overline V_i^{\sigma}\to 
\overline U_i:=\overline 
X_i-\{\{x_{i,j}\}_{j\in S_i}
\cup \{z_{i,j}\}_{j\in D_i}\}$ under the commutative finite and flat 
$\overline k$-group scheme $G_{\overline k,i}$ of rank $p$, where
$\overline X_i:=X_i\times _k \overline k$, and 
$f_i^{\sigma}$ is the image of the torsor $f_i$ above
$\overline U_i$ given by the data Deg.4, via the action of $\sigma$
on torsors (cf. [Sa], I). In other words the element $(f'_t)_t$
of $H^1_{\fppf}(U_{\overline k})_p$ determined by the datas 
Deg$^\bold {\sigma}$.4 is the transform ${(f_t)_t}^{\sigma}$
of the element $(f_t)_t$ of  $H^1_{\fppf}(U_{\overline k})_p$
determined by the datas Deg.4, under the canonical 
action of $\sigma$ on $H^1_{\fppf}(U_{\overline k})_p$ (cf. [Sa], I).

\pop {.5}
\par
\noindent
{\bf \gr Deg$^\bold {\sigma}$.5.}\rm \ For each $i\in I$ (resp. 
$i\in I^{\sp}$), are given for each $j\in S_i$ the integers 
$(m_{i,j}^{\sigma},
h_{i,j}^{\sigma})$, where $m_{i,j}^{\sigma}$
is the conductor of the torsor $f_i^{\sigma}$ at the point $x_{i,j}^{\sigma}$
of $S_i^{\sigma}$, and $h_{i,j}^{\sigma}$ its residue at this point. 
We have $m_{i,j}^{\sigma}=m_{i,j}$ and $h_{i,j}^{\sigma}=h_{i,j}$.

\pop {.5}
\par
\noindent
{\bf \gr Deg$^\bold \sigma$.6.}\rm \ For each double point 
$z_{i,j}^{\sigma}=z_{i',j'}^{\sigma}\in X_i\cap X_{i'}$ is given the pair 
$(m_{i,j}^{\sigma},h_{i,j}^{\sigma})$
(resp. $(m_{i',j'}^{\sigma},h_{i',j'}^{\sigma})$), where $m_{i,j}^{\sigma}$ 
(resp. $m_{i',j'}^{\sigma}$)
is the conductor of the torsor $f_i^{\sigma}$ (resp. $f_{i'}^{\sigma}$) 
at the point $z_{i,j}^{\sigma}$ (resp. $z_{i',j'}^{\sigma}$), 
and $h_{i,j}^{\sigma}$ (resp. $h_{i',j'}^{\sigma}$) its 
residue at this point. We have $m_{i,j}^{\sigma}=m_{i,j}$ 
and $h_{i,j}^{\sigma}=h_{i,j}$.

\pop {.5}
\par
\noindent
{\bf \gr Deg$^\bold{\sigma}$.7.}\rm \ For each double point $z_{i,j}^{\sigma}
=z_{i',j'}^{\sigma}\in 
X_i^{\sigma}\cap X_{i'}^{\sigma}$ of $\Gamma^{\sigma}$ is given the integer 
$e_{i,j}^{\sigma}=pt_{i,j}^{\sigma}$ divisible by $p$. 
For each irreducible component $X_i$ of $\Gamma$ is given the integer 
$\delta _i^{\sigma}\le v_K(p)$ which
is divisible by $p-1$, such that with the same notations as 
above $\delta _i^{\sigma}-\delta _{i'}^{\sigma}=
m_{i,j}^{\sigma}t_{i,j}^{\sigma}(p-1)$. We have $e_{i,j}^{\sigma}=e_{i,j}$,
and $\delta _i^{\sigma}=\delta _i$.

\pop {.5}
\par
\noindent
{\bf \gr Deg$^\bold \sigma$.8.}\rm \ The integer
$g_x:=(r-m-1)(p-1)/2$ is the {\bf genus} of the degeneration 
data $D_x^{\sigma}$.

\pop {.5}
\par
In conclusion, the action of the group $G_k$ on the set of isomorphism classes
of simple degeneration datas is given, and completely 
determined, through the canonical action of $G_k$ on the data Deg.4,
which is essencially given by the canonical action of $G_k$ on the group
$H^1_{\fppf}(U_{\overline k})_p$ and which was studied in [Sa].

\pop {.5}
\par
Let $\Cal X:=\Spf A$ be the formal germ of an 
$R$-curve at a smooth closed point $x$, and let
$\overline {\Cal X}:=\Spf {\overline A}$ where 
$\overline A:=A\otimes _R \overline R$. 
Let $L:=\Fr {\overline A}$ be the fractions field of $\overline A$.
In $1.2$ above we associated to a Galois cover $f:\Cal Y\to \Cal X$ of 
degree $p$, above a formal germ of an $R$-curve at a smooth point $x$, 
a simple non split (resp. split) degeneration data $D_x$. 
This indeed can be interpreted as the existence of a canonical
specialisation map:
$$\Sp:H^1_{\et}(\Spec L,\mu_p)\to \SS-Deg_p$$ 
\par
resp. 
$$\Sp:H^1_{\et}(\Spec L,\mu_p)\to \SNS-Deg_p$$ 
\par
Reciproqually, we have the following result of realisation for
degeneration datas for such covers: 

\pop {.5}
\par
\noindent
{\bf \gr 1.2.4. Theorem.}\rm \ {\sl Let $k'$ be a finite extention of $k$.
Let $\Deg(x)$ be a $k'$-simple non split 
(resp. split) degeneration data of type $(r,(G_k,m,h))$ and rank $p$. 
Then there 
exists, after eventually a finite ramified extension of $R$, a 
non-uniquely determined Galois cover $f:\Cal Y\to \Cal X$ of degree $p$, where 
$\Cal X$ is a formal germ of an $R$-curve at a smooth point $x$, and such 
that the degeneration data associated to $f$ via 1.2. is isomorphic 
to $\Deg(x)$. In other words the above specialisation maps 
$\Sp:H^1_{\et}(\Spec L,\mu_p)\to \SS-Deg_p$ and 
$\Sp:H^1_{\et}(\Spec L,\mu_p)\to \SNS-Deg_p$ are surjective.}

\pop {.5}
\par
\noindent
{\bf \gr Proof.}\ \rm One uses the formal patching results as explained in 
[Sa-1],
plus the examples given in [Sa-1] 3.1.3, 3.1.4 and 3.2.4. 
We may also assume that $k=k'$.
First, for each $i\in I^{\nsp}$ (resp. $i\in I^{\sp}$), let 
$\Cal U_i$ be a formal affine scheme with special fibre $U_i:=U_{k,i}$. 
The given torsor $f_i:V_i\to U_i$ is admissible 
(cf. [Sa], IV), hence can be lifted, after eventually a 
ramified extension of $R$, to a torsor 
$f_i':\Cal V_i\to \Cal U_i$ under a finite and flat $R$-group scheme of 
rank $p$, which is either $\mu_p$ or $\Cal H_{R,n}$ for 
$0<n\le v_K(\lambda)$ where $n:=(v_K(p)-\delta _i)/(p-1)$. Such a lifting 
is non-unique in the radicial case (cf. loc. cit). Moreover, 
for each marked point 
$x_{i,j}$ (resp. double point 
$z_{i,j}$) one can find a Galois cover $f_{i,j}:\Cal Y_{i,j}\to \Cal X_{i,j}$
of degree $p$ where  $\Cal X_{i,j}$ is a formal germ of a smooth point
$x_{i,j}$, and  $\Cal Y_{i,j}$ is smooth (resp. $f_{i,j}:\Cal Y_{i,j}\to 
\Cal Z_{i,j}$ of degree $p$,
where  $\Cal Z_{i,j}$ is a formal germ of a double point 
$z_{i,j}=z_{i'j'}$ of thikeness $e_{i,j}$, and $\Cal Y_{i,j}$ 
is semi-stable) and with 
reduction type $(G_{i,k},m_{i,j},h_{i,j})$ (resp. 
$((G_{k,i},m_{i,j},h_{i,j})$, $(G_{k,i'},m_{i',j},h_{i',j})$) on the 
boundaries) (cf. [Sa-1] 3.1.3, 3.1.4, 3.2.6). Also one can find a Galois 
cover $f_x:\Cal Y_x\to \Cal X_x$ of degree $p$, above a formal germ of a 
double point $x$, and with reduction type $(G_k,m,h)$ and 
$(G_{i_0},-m,-h)$ on the 
boundaries of the formal fibre at $x$, as 
well as a Galois cover $f'_x:\Cal Y'_x\to \Cal X'_x$ above a formal
closed disc $\Cal X'_x:=\Spf R<T>$, with degeneration data $(G_k,m,h)$
on the boundary of the formal fibre at the closed point $T=0$ (cf. [Sa-1]
2.3.1).  
Now using the formal patching result as in [Sa-1] I, as well as 
the result 3.1 in [Sa] , and after 
carefully adjusting the Galois action on the $f_{i,j}$ and the $f'_i$,
in order to obtain
Galois patching datas, one 
can patch the covers 
$f'_i$, $f_{i,j}$, $f_x$ and $f'_x$ along the points $x_{i,j}$, $z_{i,j}$ 
and $x$, 
in order to obtain a Galois cover $f':\Tilde {\Cal Y}\to \Tilde {\Cal X}$ 
of degree $p$, where $\Tilde {\Cal X}$ is a formal proper semi-stable 
$R$-curve, whose special fibre consists of the tree $\Gamma$ plus
a projective line linked to $\Gamma$ via the double point $x$. A 
$G$-equivariant contraction of the vertices $(X_i)_i$ of the 
tree $\Gamma$ in $\Tilde {\Cal X}_k$, will yield to a Galois cover 
$f:\Cal Y'\to \Cal X'$ of degree $p$,
where $\Cal X'$ is a formal proper and smooth projective line, 
with a marked point $x$ on $\Cal X'_k$. A formal localisation now
at the point $x$ will give the desired cover $f:\Cal Y_x\to \Cal X_x
:=\Spf \hat \Cal O_{\Cal X',x}$, and by construction the simple 
degeneration data associated to $f$ via 1.2 is isomorphic to $D_x$.

\pop {.5}
\par
\noindent
{\bf \gr 1.2.5. Example.}\rm\  This in particular is the realisation given 
by the above theorem 1.2.4 of the degeneration datas which arises in the 
example 1.2.1. Consider the non split 
simple degeneration data of rank $p$ which consists of a tree with one 
vertex and no edges. Hence a projective line $X_1$ with one 
$\overline k$-marked point $x=x_1$, and a given \'etale torsor 
$f_1:V_1\to U_1:=X_1-\{x_1\}$ with conductor $m'-m+1$ at the point $x_1$
where $m'>m$ are positif integers. Then by 1.2.4 one can construct 
(after eventually a finite extension of $R$) a Galois 
cover $f:\Cal Y\to \Cal X$ of degree $p$ above a formal germ $\Cal X$ 
at a smooth point $x$, such that the singularity of $\Cal Y$ at the closed 
point $y$ is unibranche and the geometric genus of the point $y$ equals 
$(m'-m)(p-1)/2$, moreover the semi-stable reduction of $\Cal Y$ consists of
one irreducible component $Y_1$ of genus $(m'-m)(p-1)/2$, 
which is the projective completion of the above affine curve $V_1$, 
and which is linked to the point $y$ by a double point.

\pop {.5}
\par
\noindent
{\bf \gr 1.2.6.}\rm \ Now we will study the action of $G_K$ on Galois 
covers of degree $p$ above formal germs of smooth 
curve at closed points. This action as we will see extends in a 
functorial way to an action of $G_K$ on the corresponding 
semi-stable models. More precisely, let $\Cal X:=\Spf A$ be the 
formal germ of an $R$-curve at a smooth closed point $x$, and let
$\overline {\Cal X}:=\Spf {\overline A}$ where $\overline A:=A
\otimes _R \overline R$. 
Let $L:=\Fr {\overline A}$ be the fractions field of $\overline A$. 
Then the Galois group
$G_K$ acts in a natural way on cyclic extentions of degree $p$ of $L$, 
and we have a canonical homomorphism: 
$$(1)\ G_K\to \Aut H^1_{\et}(\Spec L,\mu_p)$$

\par
The above action $(1)$ is equivalent to the action of $G_K$ on 
normal Galois covers of degree $p$ above $\overline {\Cal X}$. 
More precisely to each such a cover 
$f:{\Cal Y}\to \overline {\Cal X}$ with 
${\Cal Y}$ local and normal, which corresponds to the extension of 
functions fields: $\Spec L'\to \Spec L$, 
and an element $\sigma \in G_K$, one associates the Galois cover 
$f^{\sigma}:{\Cal Y}^{\sigma}\to \overline {\Cal X}$ 
of degree $p$, with 
${\Cal Y}^{\sigma}$ local and normal, 
which corresponds to the extension of 
functions fields: ${\Spec L'}^{\sigma}\to \Spec L$.
This action of $G_K$ on covers extends to an action of $G_K$ on 
corresponding minimal semi-stable models (because of the minimality 
condition). Hier by a minimal semi-stable model we mean either a split 
or a non split one. More precisely, let  
$\Tilde f:\Tilde {\Cal Y}\to {\Cal Y}$ (resp. 
$\Tilde f':\Tilde {\Cal Y}'\to 
{\Cal Y}^{\sigma}$) be a 
minimal semi-stable model of ${\Cal Y}$ (resp. of
${\Cal Y}^{\sigma}$), then the element $\sigma\in G_K$ maps the 
semi-stable model $\Tilde f$ to the semi-stable model $\Tilde f'$,
in the sens that we have a commutative diagram:

$$
\CD
\Tilde {\Cal Y}    @>\Tilde f>>    {\Cal Y}    @> f>> \overline {\Cal X} \\
   @V{\sigma}VV                  @V{\sigma}VV           @VVV\\
\Tilde {\Cal Y}'   @>{\Tilde f}^{\sigma}>>  
{\Cal  Y}^{\sigma} 
@>f^{\sigma}>>  {\overline {\Cal X}}
\endCD
$$

where the horizontal arrows are $\overline {\Cal X}
$-automorphisms. We first show how the above action of $G_K$ induces an 
action of $G_K$ on simple degeneration datas via its canonical quotient
$G_k$. Our main result then is that this action coincides with the canonical 
action of $G_k$ on degeneration datas as defined in 1.2.2, and that the 
specialisation maps $\Sp:H^1_{\et}(\Spec L,\mu_p)\to \SS-Deg_p$ and 
$\Sp:H^1_{\et}(\Spec L,\mu_p)\to \SNS-Deg_p$ are $G_K$-equivariant.

\pop {.5}
\par
First we may assume that the Galois cover $f:\Cal Y\to \Cal X':=
\Cal X\times _RR'$,
as well as its Galois conjugate via $\sigma$, are defined over a finite 
extension $R'$ of $R$. Let $\pi '$ be a uniformiser of
$R'$, and let $A':=A\otimes _R R'$. Let $\wp:=(\pi')$ be the ideal of 
$A'$ generated by $\pi'$, and let $\hat A'_{\wp}$ be the 
completion of the localisation of $A'$ at $\wp$. Let 
$\Cal X'_{\eta}:=\Spf \hat A'_{\wp}$ be the boundary of $\Cal X'$, and let 
$\Cal X'_{\eta}\to \Cal X'$ be the canonical morphism. 
Consider the following cartesian diagram:

$$
\CD
\Cal Y    @> f>> \Cal X' \\
   @VVV                  @VVV   \\
{\Cal  Y}_{\eta} @>f_{\eta}>>  {\Cal X}'_{\eta}
\endCD
$$

Then $f_{\eta}:\Cal Y_{\eta}\to \Cal X'_{\eta}$ is a torsor under a 
commutative finite and flat $R'$-group scheme $G_{R'}$ 
of rank $p$. Let $(G_{k'},m,h)$ be the degeneration type of $f_{\eta}$
where $k'$ is the residue field of $R'$.
We have also the following cartesian diagram:

$$
\CD
\Cal Y  ^{\sigma}  @> f^{\sigma}>> \Cal X' \\
   @VVV                  @VVV   \\
{\Cal  Y_{\eta}^{\sigma}} @>f_{\eta}^{\sigma}>>  {\Cal X}'_{\eta}
\endCD
$$

\par
Moreover the torsor $f_{\eta}^{\sigma}$ is the Galois
conjugate of the torsor $f_{\eta}$ by $\sigma$ as one sees easily.
In particular both of these torsors have the same reduction type.
\pop {.5}
\par
The Galois group $G\simeq \Bbb Z/p\Bbb Z$ of the cover $f$ (resp. 
$f^{\sigma}$) acts in a canonical way on $\Tilde {\Cal Y}$ (resp. 
$\Tilde {\Cal Y}':=\Tilde {\Cal Y}^{\sigma}$). Let 
$\Tilde {\Cal X}:=\Tilde {\Cal Y}/G$ which also 
equals $\Tilde {\Cal Y}'/G$. We have a commutative diagram:

$$
\CD
\Tilde {\Cal Y}    @>>> \Tilde {\Cal X} \\
   @V\sigma VV                  @VVV   \\
\Tilde {\Cal Y}^{\sigma} @>>>  \Tilde {\Cal X}
\endCD
$$

Where the horizontal arrows are Galois covers with group $G$, and
$\sigma$ is an 
\newline
$\Tilde {\Cal X}$-automorphism. Let
$\Deg(x_1)$ (resp. $\Deg(x_2)$) be the simple degeneration data
associated to the cover $f$ (resp. $f^{\sigma}$). Let $\overline
\sigma$ be the image of $\sigma$ in $G_k$ via the canonical surjective 
homomorphism $G_K\to G_k$. We will prove that $\Deg (x_2)$ is the transform
$\Deg(x_1)^{\overline {\sigma}}$ of $\Deg(x,1)$ by $\overline {\sigma}$, 
via the canonical action of $G_k$ on degeneration datas as explained in
[Sa], I. The only think we have to check to prove 
this is that this is true when 
considering the degenertaion data Deg.4. Let $X_i$ be an 
irreducible component of $\Tilde {\Cal X}_{k'}:=
\Tilde {\Cal X}\times _R k'$. Let $D_i$ be the set
of double points of $\Tilde {\Cal X}_{k'}$ supported by $X_i$, and let
$S_i$ be the set of smooth points of $X_i$ in which specialise some branched 
points of $f$. Let $U_i:=X_i-\{S_i\cup D_i\}$, and let 
$\Cal U_i$ be a formal open
subset of $\Tilde {\Cal X}$ with special fibre $U_i$. Let 
$f_i:\Cal V_i\to \Cal U_i$ (resp. $f_i':\Cal V_i'\to \Cal U_i$)
be the torsor above $\Cal U_i$ induced by $f$ (resp. by $f^{\sigma}$).
Then $f_i'$ is the transform $f_i^{\sigma}$ by $\sigma$ of the torsor
$f_i$, via the canonical action of $G_K$ on these torsors. Let
$f_{i,k'}:\Cal V_{i,k'}\to \Cal U_{i,k'}$ 
(resp. $f_{i,k'}':\Cal V_{i,k'}'\to \Cal U_{i,k'}$) be the special fibre
of the torsor $f_i$ (resp. $f_i':=f_i^{\sigma}$). Let $\overline {\sigma}$ 
be the image of $\sigma$ in $G_k$ via the canonical homomorphism
$G_K\to G_k$. Then $f_{i,k'}'$ is nothing else but the transform 
${f_{i,k'}}^{\overline {\sigma}}$ of $f_{i,k'}$ via ${\overline {\sigma}}$
(cf. [Sa], I). So indeed we have proven the following:

\pop {.5}
\par
\noindent
{\bf \gr 1.2.7. Theorem.}\rm \ {\sl Let $G_K$ acts on the sets 
$\SS-Deg_p$ (resp. $\SNS-Deg_p$) of isomorphism classes of 
simple split (resp. simple non split) degeneration data of rank $p$, 
via its canonical
quotient $G_k$, and through the natural action of 
$G_k$ on these sets as defined in 1.2.3. Then the surjective specialisation 
maps $\Sp:H^1_{\et}(\Spec L,\mu_p)\to \SS-Deg_p$ and 
$\Sp:H^1_{\et}(\Spec L,\mu_p)\to \SNS-Deg_p$ are $G_K$-equivariant.}

\pop {.5}
\par
\noindent
{\bf \gr 1.2.8. Remark.}\rm \ It is easy to construct examples of covers 
$f:\Cal Y\to \Cal X$ as in [Sa-1] 3.1.4, where the special fibre $\Cal Y_k$ is 
singular and unibranche at the closed point $y$ of $\Cal Y$, and such that
the configuration of 
the special fibre of a semi-stable model $\Tilde {\Cal Y}$ of $\Cal Y$ is 
not a tree-like. This indeed answers a question raised by Lorenzini
whether such a situation can occur in the inequal characteristic case.
More precisely, consider the non split simple degeneration data $\Deg (x)$
of type $(2,(\mu_p,-1,0))$ which consists of a graph $\Gamma$ with two
vertices $X_1$ and $X_2$ linked by a unique edge $y$, with given 
$\overline k$-marked points $x=x_1$ on $X_1$ and  $\overline k$-marked 
point $x_2$ on $X_2$, and given \'etale torsors of rank $p$ :
$f_1:V_1\to U_1:=X_1-\{x_1\}$ with conductor $m_1=1$ at $x=x_1$ and  
$f_2:V_2\to U_2:=X_2-\{x_2\}$ with conductor $m_2=1$ at $x_2$. Then it 
follows from 1.2.4 that there exists, after eventually a finite extension 
of $R$, a Galois cover $f:\Cal Y\to \Cal X$ of degree $p$ above the formal 
germ $\Cal X\simeq \Spf R[[T]]$ at the smooth $R$-point $x$, such that the 
$\overline k$-simple non split degeneration data associated to the above 
cover $f$ is the above 
given one. Moreover by construction the singularity of the closed point $y$ 
of $\Cal Y$ is unibranche, and the configuration of the (non-split) 
semi-stable reduction of $\Cal Y$ consists of two projective lines which 
meet at $p$-double points (the above cover will be \'etale in reduction 
above the double point $y$), in particular one has $p-1$ cycles in this 
configuration.

\vskip.4cm
\epsfysize=7.4cm
\centerline{\epsfbox{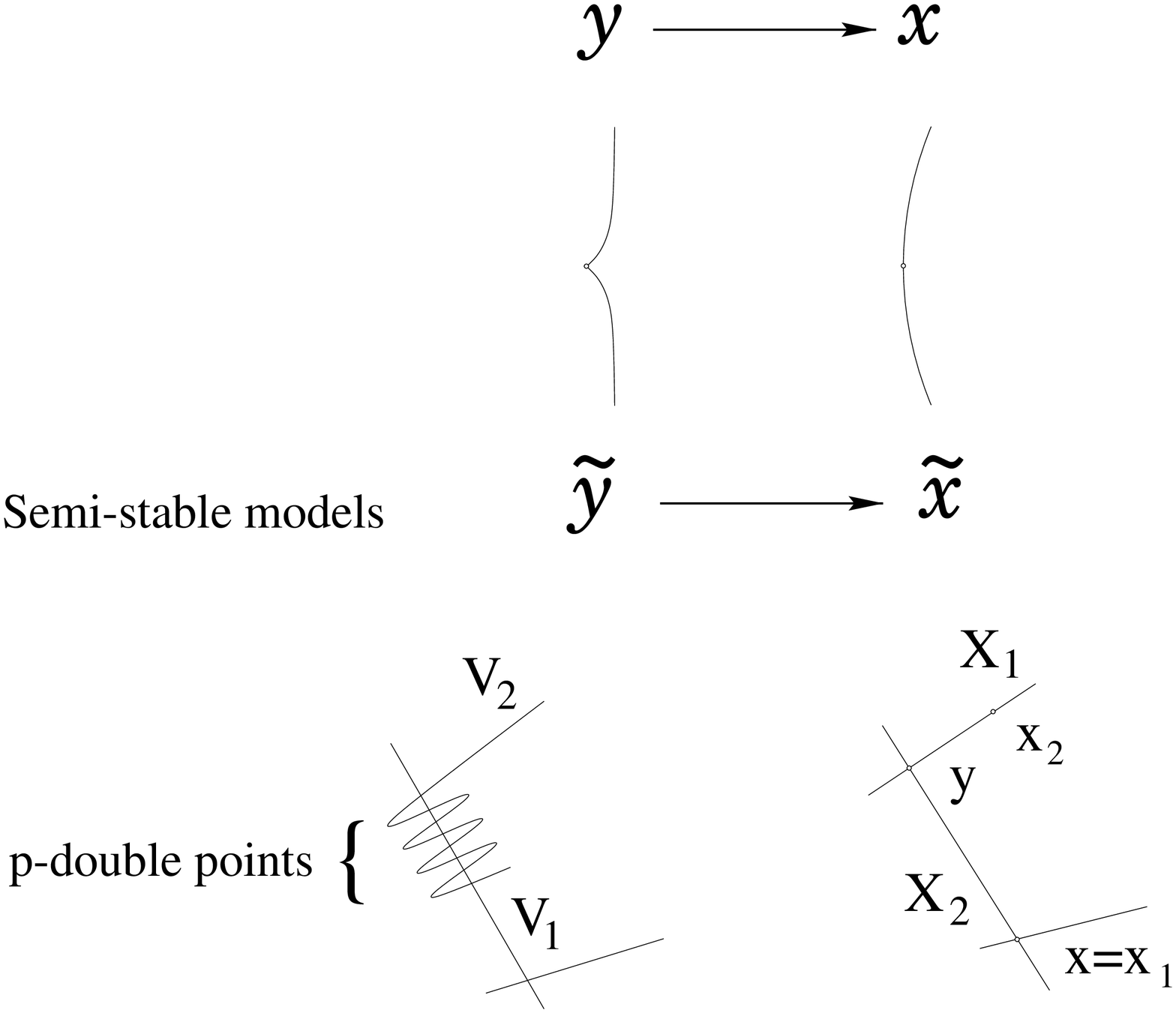}}

\pop {1}
\par
\noindent
{\bf \gr 1.3.}\rm \ We use the same notations as in 1.1. 
We consider now the case where $x$ is an {\bf ordinary 
double} point. Let $f:\Cal Y\to \Cal X$ be a Galois cover of degree $p$,
where $\Cal X$ is the formal germ of an $R$-curve at an ordinary double 
point $x$. In a similar way as in 1.2, and using the results of [Sa] and 
[Sa-1], one can associate a {\it double degeneration} data to $f$ 
whose definition we summarise in the following:

\pop {.5}
\par
\noindent
{\bf \gr 1.3.1. Definition.}\rm \  Let $k'$ be an algebraic extension
of $k$. A $k'$-{\bf double non split (resp. split) degeneration data}
$\Deg (x)$ of type $(r,(G_{k',1},m_1,h_1),(G_{k',2},m_2,h_2))$ and rank $p$
consists of the following:

\pop {.5}
\par
\noindent
{\bf \gr Deg.1.}\rm \ For $i\in \{1,2\}$, $G_{k',i}$ is a finite and flat 
$k'$-group scheme of rank $p$, 
$r\ge 0$ is an integer, $m_i$ is an integer prime to $p$ such that 
$r-m_1-m_2\ge 0$, and $h_i\in \Bbb F_p$ equals $0$ unless 
$G_{k',i}=\mu _p$ and $m_i=0$.

\pop {.5}
\par
\noindent
{\bf \gr Deg.2.}\rm\  $\Gamma:=X_{k'}$ is an oriented tree of projective lines 
with vertices $\Vert (\Gamma):=(X_i)_{i\in I}$ 
(resp. $\Vert (\Gamma):=(X_i)_{i\in I^{\sp}}$), and 
$\Tilde {\Gamma}=\cup _{i\in \Tilde I} X_{i}$, where 
$\Tilde I\subset I$, is a geodesic in $\Gamma$ linking two given vertices
$X_{i_1}$ and $X_{i_2}$, with a given marked smooth $k'$-rational 
point $x_1\in X_{i_1}$ (resp. $x_2\in X_{i_2}$). We denote by 
$\{z_{i,j}\}_{j\in D_i}$ the set of double points, or edges, of $\Gamma$ 
which are supported by $X_i$.

\pop {.5}
\par
\noindent
{\bf \gr Deg.3.}\rm \ For each vertex $X_i$ of $\Gamma$ is given a set, 
may be empty, of smooth marked $k'$-rational points 
$\{x_{i,j}\}_{j\in S_i}$. 

\pop {.5}
\par
\noindent
{\bf \gr Deg.4.}\rm \  For each $i\in I$ (resp. $i\in I^{\sp}$), is given a 
torsor $f_{i,k'}:V_{i,k'}\to U_{i,k'}:= 
X_{i,k'}-\{\{x_{i,j}\}_{j\in S_i}
\cup \{z_{i,j}\}_{j\in D_i}\}$ under a commutative finite and flat 
$k'$-group scheme $G_{k',i}$ of rank $p$, where
$X_{i,k'}:=X_i\times _{\Spec k}\Spec {k'}$, 
with $V_{i,k'}$ {\bf smooth}. Let $U_{k'}$ be the open subscheme
of $X_{k'}$ obtained by deleting the double and marked points
of $X_{k'}$. Then the above data determines an element of
$H^1_{\fppf}(U_{k'})_p$. Moreover if $G_{k',i}$ is 
radicial, and if $\omega _i$ is the associated differential form, then 
the set of zeros and poles of $w_i$ is necessarily contained in  
$\{x_{i,j}\}_{j\in S_i}\cup \{z_{i,j}\}_{j\in D_i}$, as
$V_{i,k'}$ is smooth. Also if $i\in I^{\sp}$, and $S_i$ is non empty, we 
assume that $G_{i,k'}=\mu_p$.

\pop {.5}
\par
\noindent
{\bf \gr Deg.5.}\rm \ For each $i\in I$ (resp. $i\in I^{\sp}$), are given for 
each $j\in S_i$ integers $(m_{i,j},h_{i,j})$, where $m_{i,j}$
is the conductor of the torsor $f_i$ at the point $x_{i,j}$, 
and $h_{i,j}$ its residue at this point. Let $r_{i,j}:=m_{i,j}+1$.
If $i\in I^{\sp}$ then $r_{i,j}=1$.

\pop {.5}
\par
\noindent
{\bf \gr Deg.6.}\rm \  For each double point 
$z_{i,j}=z_{i',j'}\in X_i\cap X_{i'}$ is given the pair $(m_{i,j},h_{i,j})$
(resp. $(m_{i',j'},h_{i',j'})$), where $m_{i,j}$ (resp. $m_{i',j'}$)
is the conductor of the torsor $f_{i,k'}$ (resp. $f_{i',k'}$) at the point 
$z_{i,j}$ (resp. $z_{i',j'}$), and $h_{i,j}$ (resp. $h_{i',j'}$) its 
residue at this point. These datas must satisfy $m_{i,j}+m_{i',j'}=0$, 
and $h_{i,j}+h_{i',j'}=0$. In other words the element
of $H^1_{\fppf}(U_{k'})_p$ determined by the data 
Deg.4 must belong
to the subgroup $H^1_{\fppf}(U_k)_p^{\kum}$. In particular 
for the double point 
$x_{i_1,j_1}:=x_1$ (resp. $x_{i_2,j_2}:=x_2$)
is given the data $(m_{i_1,j_1},h_{i_1,j_1})$ (resp.
$(m_{i_2,j_2},h_{i_2,j_2}$), such that $m_{i_1,j_1}+m_1=0$ 
(resp. $m_{i_2,j_2}+m_2=0$), and 
$h_{i_1,j_1}+h_1=0$ (resp. $h_{i_2,j_2}+h_2=0$).

\pop {.5}
\par
\noindent
{\bf \gr Deg.7.}\rm \ For each irreducible component
$X_i$ of $\Gamma$ is given an integer $\delta _i\le v_K(p)$ which
is divisible by $p-1$. For each double point $z_{i,j}=z_{i',j'}\in 
X_i\cap X_{i'}$ of $\Gamma$ is given an integer 
$e_{i,j}=pt_{i,j}$ divisible by $p$, such that with the same notations as 
above we have $\delta _i-\delta _{i'}=m_{i,j}t_{i,j}(p-1)$.

\pop {.5}
\par
\noindent
{\bf \gr Deg.8.}\rm \ Let $I_{\et}$ denote the subset of $I$ 
(resp. of $I^{\sp}$) consisting of those $i$ for which
$G_{k',i}$ is \'etale, then : $(r-m-1)(p-1)/2
=\sum _{i\in I_{\et}}(-2+\sum _{j\in S_i}
(m_{i,j}+1)+\sum _{j\in D_i}(m_{i,j}+1))(p-1)/2$. The integer 
$g_x:=(r-m_1-m_2)(p-1)/2$ is called the {\bf genus} of the 
degeneration data $\Deg (x)$.

\pop {.5}
\par
There is a natural notion of isomorphism of double split (resp. non split) 
degeneration datas of a given type over a given algebraic extension of
$k$. We will denote by $\bold {\DS-Deg_p}$ (resp. 
$\bold {\DNS-Deg_p}$) the set of 
isomorphism classes of $\overline k$-double split (resp. 
$\overline k$-double non split) degeneration data of rank $p$. 
Also as in 1.2.3 one shows that there exists a canonical action of 
$G_k$ on these sets, via the canonical action of
$G_k$ on the degeneration data Deg.4, and these sets are 
$G_k$-sets, hence also canonically $G_K$-sets.
The following theorem is proven  in the same way as in 1.2.4 and 1.2.6.

\pop {.5}
\par
\noindent
{\bf \gr 1.3.2. Theorem.}\rm \ {\sl Let $\Cal X:=\Spf A$ be the formal 
germ of an $R$-curve at an ordinary double point. Let $\overline A:=
A\otimes _R \overline R$, and let $L:=\Fr (\overline A)$. Then the canonical
specialisation maps: $\Sp:H^1_{\et}(\Spec L,\mu_p)\to \DS-Deg_p$ and 
$\Sp:H^1_{\et}(\Spec L,\mu_p)\to \DNS-Deg_p$ are surjective and
$G_K$-equivariant.

\pop {.5}
\par
\noindent
{\bf \gr 1.3.3. remark.}\rm \ The surjectivity result in 1.2.4 and 1.2.7
was shown in [He] in the case of degeneration datas of genus $0$ 
in the framework of automorphisms of order $p$ of open discs and annuli. Our
method based on the techniques developed in [Sa] and [Sa-1] avoid the use of 
this language and uses only kummer theory.

\pop {1}
\par
\noindent
{\bf \gr II. Semi-stable reduction of Galois covers of degree $p$
obove proper curves.}
\rm

\pop {.5}
\par
\noindent
{\bf \gr 2.0.}\rm\ In this section we will use the same notations as in 
1.0. We will consider a {\bf proper} and {\bf semi-stable} $R$-curve $X$, 
a Galois cover $f:Y\to X$ with group $\Bbb Z/p\Bbb Z$, and with 
$Y$ normal. We will study the semi-stable reduction of $Y$ as well as 
the Galois action on these datas.

\pop {.5}
\par
\noindent
{\bf \gr 2.1.}\rm\ Let $X$ be a proper and semi-stable $R$-curve,
with smooth and connected generic fibre $X_K:=X\times _RK$. We assume that 
the double points of $X_k$ are $k$-rationals. Let $f:Y\to X$ be a Galois 
cover with group $G\simeq \Bbb Z/p\Bbb Z$, and with 
$Y$ normal. It follows from the
theorem of semi-stable reduction for curves (cf. [De-Mu]) that, after 
eventually a finite 
extension $R'$ of $R$, with fractions field $K'$ and residue field $k'$, 
the curve $Y$ has a semi-stable reduction. More 
precisely, there exists a birational and proper morphism 
$\Tilde f :\Tilde {Y}\to Y':=Y\times _R R'$ such that $\Tilde {Y}_{K'}\simeq 
 Y'_{K'}$, and 
the following conditions hold:
\par
\noindent
(i) The special fibre $\Tilde {Y}_{k'}$ of $\Tilde {Y}$ is reduced.
\par
\noindent
(ii) $\Tilde {Y}_{k'}$ has only ordinary double points as singularities.

\par
Moreover there exists such a semi-stable model $\Tilde {Y}$ which is
minimal. In particular the action of $G$ on $Y'$ extends to an action
on $\Tilde {Y}$. Let $\Tilde {X}$ be the quotient of 
$\Tilde {Y}$ by $G$ which is a semi-stable model of $X':=X\times _R R'$. 
In the following we will choose $\Tilde {Y}$ and $\Tilde {X}$ as above 
such that the
set of points $B_K:=\{x_{i,K}\}_{i=1}^r$ consisting of the branch locus of the 
morphism $f_K:Y_K\to X_K$ are rational. Moreover we choose 
such $\Tilde {X}$ 
and $\Tilde {Y}$ which are minimal for these properties. One has the 
following commutative diagram:

$$
\CD
\Tilde Y @>g>>   \Tilde X \\
    @V\Tilde fVV           @V\Tilde g VV   \\
Y' @>f'>>      X'
\endCD
$$

\par 
Let $(X_i)_{i\in I}$ be the irreducible components of the special fibre 
$X'_{k'}$ of $X'$. For each $i\in I$, let $(x_{i,j})_{j\in S_i}$ be the set
(may be empty) of those smooth points of $X_i$ in which specialise some 
points of $B_K$, say in $x_{i,j}$ specialise $r_{i,j}$ points of $B_K$, and
let $(z_t)_{t\in J}$ be the set of double points of $X_k$. For each 
$t\in J$, let $r_t$ be the number of points of $B_K$ which specialise in
$z_t$, we have $r=\sum _{t\in J}r_t+\sum _{i\in I}\sum _{j\in S_j}r_{i,j}$.
Let $U:=X'-\{D\cup S\}$, where $S:=\cup _{i\in I}\{x_{i,j}\}_{j\in S_i}$
and $D=\{z_t\}_{t\in J}$. We will denote by $\{U_i\}_{i\in I}$ the set of 
irreducible components of $U$.

\pop {.5}
\par
\noindent
{\bf \gr 2.2.}\rm\ Below we will use the notations of 2.0 and 2.1 and will
exhibit the degeneration datas associated to the semi-stable reduction of $Y$.
Let $R'$ be a finite extension of $R$ such that the conditions of 2.1 hold.
The following datas are canonically associated to the Galois cover $f:Y\to X$:

\pop {.5}
\par
\noindent
{\bf \gr Deg.1.}\rm \ For each irreducible component $X_i$ of $X'_{k'}$,
let $\eta_i$ be the generic point of $X_i$, and let $X'_{\eta_i}$ be the 
spectrum of the completion of the localisation of $X'$ at $\eta_i$. We have 
a canonical morphism : $X'_{\eta_i}\to X'$. Consider the following cartesian 
diagram:

$$
\CD
Y_{\eta_i}  @>f_{\eta_i}>>   X'_{\eta_i}  \\
    @VVV             @VVV   \\
Y'      @>f'>>        X'
\endCD
$$

Then either $f_{\eta_i}:Y_{\eta_i}\to X'_{\eta_i}$ is completely split, or
is a torsor under a finite and flat commutative $R'$-group scheme $G_i$ 
of rank $p$.

\pop {.5}
\par
\noindent
{\bf \gr Deg.2.}\rm \ For each irreducible component $X_i$ of $X'_k$, let
$U_i:=X_i-\{\{x_{i,j}\}_{j\in S_i}\cup \{z_t\}_{t\in J_i}\}$, where 
$J_i\subset J$ 
denotes the index subset indexing those double points of $X'_{k'}$ 
which are supported by
$X_i$. Then $V_i:=f^{-1}(U_i)\to U_i$ is an admissible torsor under the 
finite and flat group scheme $G_{k',i}:={G_i}\times_{R'}{k'}$. Moreover, if 
$G_{k',i}$ is \'etale then $V_i$ is smooth, and if $G_{k',i}$ is radicial 
then the only singularities of $V_i$ lie above the zeros 
$\{y_l\}_{l\in Z_i}$ of the differential form $\omega_i$ associated to 
the above torsor (cf. [Sa], I).

\pop {.5}
\par
\noindent
{\bf \gr Deg.3.}\rm \ For each double point $z_j\in X_i\cup X_{i'}$, let
$\Cal X_j$ be the completion of the localisation of $X'$ at $z_j$, and let 
$\Cal X_j\to X'$ be the canonical morphism. Consider the cartesian diagram:

$$
\CD
\Cal Y_{j}  @>f_j>>   \Cal X_{j}  \\
    @VVV             @VVV   \\
Y'      @>f'>>        X'
\endCD
$$

Then either the cover $f_j:\Cal Y_j\to \Cal X_j$ is split, or it is a Galois 
cover of degree $p$ with $\Cal Y_j$ connected, in which case is associated 
to $f_j$ a double degeneration data $\Deg (z_j)$. Let $m_{j,1}$ and $h_{j,1}$
(resp. $m_{j,2}$ and $h_{j,2}$) be the conductor and the residue at the 
double point $z_j$ associated to the torsor $V_i:=f^{-1}(U_i)\to U_i$
(resp. $V_{i'}:=f^{-1}(U_{i'})\to U_{i'}$). Then the above degeneration 
data $\Deg (z_i)$ is of type $(r_j,(G_{k',i},m_{j,1},h_{j,1}),
(G_{k',i'},m_{j,2},h_{j,2}))$, and genus 
$(r_j-m_{j,1}-m_{j,2})(p-1)/2$.

\pop {.5}
\par
\noindent
{\bf \gr Deg.4.}\rm \  For each irreducible component $X_i$ of 
$X'_{k'}$, and a smooth point $x_{i,j}$, $j\in S_i$, let $\Cal X_{i,j}$ be 
the completion of the localisation of $X'$ at $x_{i,j}$, and let 
$\Cal X_{i,j}\to X'$ be the canonical morphism. Consider the cartesian 
diagram:

$$
\CD
\Cal Y_{i,j}  @>f_{i,j}>>   \Cal X_{i,j}  \\
    @VVV             @VVV   \\
Y'      @>f'>>        X'
\endCD
$$
The cover $f_{i,j}:\Cal Y_{i,j}\to \Cal X_{i,j}$ is 
necessarily a Galois cover of 
degree $p$ with $\Cal Y_{i,j}$ connected, and to $f_{i,j}$ is associated 
a simple degeneration data $\Deg (x_{i,j})$. Let $m_{i,j}$ and 
$h_{i,j}$ be the conductor and the residue at the 
smooth point $x_{i,j}$ associated to the torsor $V_i:=f^{-1}(U_i)\to U_i$.
Then the above degeneration data $\Deg (x_{i,j})$ is of type 
$(r_{i,j},(G_{k',i},m_{i,j},h_{i,j}))$, and genus 
$(r_{i,j}-m_{i,j}-1)(p-1)/2$.

\pop {.5}
\par
\noindent
{\bf \gr Deg.5.}\rm \ For each $i\in I$ such that the group scheme 
$G_{k',i}$ is radicial, let $\{y_l\}_{l\in Z_i}$ be the zeros of the 
differential form $\omega_i$ associated to the above torsor 
$V_i:=f^{-1}(U_i)\to U_i$, and which we may assume to be rational over $k'$. 
For each $l\in Z_i$,  let $\Cal Z_{l}$ be 
the completion of the localisation of $X'$ at $y_{l}$, and let 
$\Cal Z_{l}\to X'$ be the canonical morphism. Consider the cartesian diagram:

$$
\CD
\Cal Z'_{l}  @>f_{l}>>   \Cal Z_{l}  \\
    @VVV             @VVV   \\
Y'      @>f'>>        X'
\endCD
$$
The cover $\Cal Z'_{l}\to \Cal Z_{l}$ is a Galois cover of 
degree $p$ with $\Cal Z'_{l}$ connected, and to $f_{l}$ is associated 
a simple degeneration data $\Deg (y_{l})$. Let $m_{l}$ and 
$h_{l}$ be the conductor and the residue at the 
smooth point $y_l$ associated to the torsor $V_i:=f^{-1}(U_i)\to U_i$.
Then the above degeneration data $\Deg (y_l)$ is of type 
$(0,(G_{k',i},m_{l},h_{l}))$, and genus 
$(-m_{l}-1)(p-1)/2$.

\pop {.5}
\par
The above considerations lead to the following definition of degeneration 
datas associated to a Galois cover of degree $p$ above a proper 
and semi-stable $R$-curve.

\pop {.5}
\par
\noindent
{\bf \gr 2.2.1. Definition.}\rm \ Let $X_k$ be a proper and semi-stable 
$k$-curve, and let $\{X_i\}_{i\in I}$ be the irreducible components 
of $X_k$. We assume that the double points
$\{z_t\}_{t\in J}$ of $X_k$ are $k$-rational. Let $k'$ be an algebraic 
extension $k$. A $k'$-{\bf degeneration data of rank $p$} associated to $X_k$, 
consists of the following datas:

\pop {.5}
\par
\noindent
{\bf \gr Deg.1.}\rm \ For each irreducible component $X_i$ of $X_k$
is given a set of smooth $k'$-rational points $\{x_{i,j}\}_{j\in S_i}$ 
of $X_k$ which are supported by $X_i$.
\pop {.5}
\par
\noindent
{\bf \gr Deg.2.}\rm \ For each irreducible component $X_i$ of $X_k$, let
$U_i:=X_i-\{\{x_{i,j}\}_{j\in S_i}\cup \{z_t\}_{t\in J_i}\}$, where 
$J_i\subset J$ 
denotes the index subset indexing those double points of $X_k$ 
which are supported by $X_i$. Then we assume given an admissible torsor
$f_i:V_i\to U'_i:=U_i\times _k k'$ (cf. [Sa], 4, for the definition 
of admissibility) under a finite and flat $k'$-group scheme 
$G_{k',i}$, and we allow the torsor $f_i$ to be trivial. If $G_{k',i}$ is 
radicial then we assume that the zeros 
$\{y_l\}_{l\in Z_i}$ of the differential form $\omega_i$ associated to 
the above torsor $f_i$ are $k'$-rational. 
\pop {.5}
\par
\noindent
{\bf \gr Deg.3.}\rm \  For each double point $z_j\in X_i\cap X_{i'}$ of 
$X_k$, let $m_{j,1}$ and $h_{j,1}$ (resp. $m_{j,2}$ and $h_{j,2}$) be 
the conductor and 
the residue at the double point $z_j$ associated to the torsor 
$V_i:=f^{-1}(U_i)\to U_i$ (resp. $V_{i'}:=f^{-1}(U_{i'})\to U_{i'}$). 
Let $r_j$ be an integer such that $r_j-m_{j,1}-m_{j,2}\ge 0$. Then we assume 
given a $k'$-double degeneration data
$\Deg (z_j)$ of type $(r_j,(G_{k',i},m_{j,1},h_{j,1}),
(G_{k',i'},m_{j,2},h_{j,2}))$, and genus $(r_j-m_{j,1}-m_{j,2})(p-1)/2$.
\pop {.5}
\par
\noindent
{\bf \gr Deg.4.}\rm \  For each irreducible component $X_i$ of 
$X_{k'}$, and a smooth point $x_{i,j}$, $j\in S_i$, let $m_{i,j}$ and 
$h_{i,j}$ be the conductor and the residue at the 
smooth point $x_{i,j}$ associated to the torsor $V_i:=f^{-1}(U_i)\to U_i$.
Let $r_{i,j}$ be an integer such that $r_{i,j}-m_{i,j}-1\ge 0$. Then we 
assume given
a $k'$-simple degeneration data $\Deg (x_{i,j})$ of type 
$(r_{i,j},(G_{k',i},m_{i,j},h_{i,j}))$, and genus $(r_{i,j}-m_{i,j}-1)(p-1)/2$.
\pop {.5}
\par
\noindent
{\bf \gr Deg.5.}\rm \ For each $i\in I$ such that the group scheme 
$G_{k',i}$ is radicial, let $\{y_l\}_{l\in Z_i}$ be the zeros of the 
differential form $\omega_i$ associated to the above torsor 
$V_i:=f^{-1}(U_i)\to U_i$, and which we may assume to be rational over $k'$. 
For each $l\in Z_i$, let $m_{l}$ and 
$h_{l}$ be the conductor and the residue at the 
smooth point $y_l$ associated to the torsor $V_i:=f^{-1}(U_i)\to U_i$.
Then we assume given a degeneration data $\Deg (y_l)$ of type 
$(0,(G_{k',i},m_{l},h_{l}))$, and genus $(-m_{l}-1)(p-1)/2$.
\pop {.5}
\par
\noindent
{\bf \gr 2.2.2. Remark.}\rm \ The admissibility condition  
on the torsors $\{f_i\}_{i\in I}$ which is part of the above definition of
degeneration datas is satisfied when the $U_i$ are opens 
of the projective line, e.g. if $X_k$ is a totally degenerate Mumford curve, 
and in the case where $U_i$ is proper.

\pop {.5}
\par
\noindent
{\bf \gr 2.2.3. Example.} \rm In this example we consider a semi-stable 
and proper $R$-curve $X$ with smooth and geometrically connected generic 
fibre $X_K$, and whose special fibre $X_k$ consists of two irreducible 
smooth components $X_1$ of genus $g_1\ge 2$ and $X_2$ of genus $g_2\ge 2$
which intersect at the double point $x$. We suppose moreover that both 
$X_1$ and $X_2$ are {\bf generic} (for this one has to choose $R$ 
to be rather large). 

\pop {1}
\epsfysize=3cm
\centerline{\epsfbox{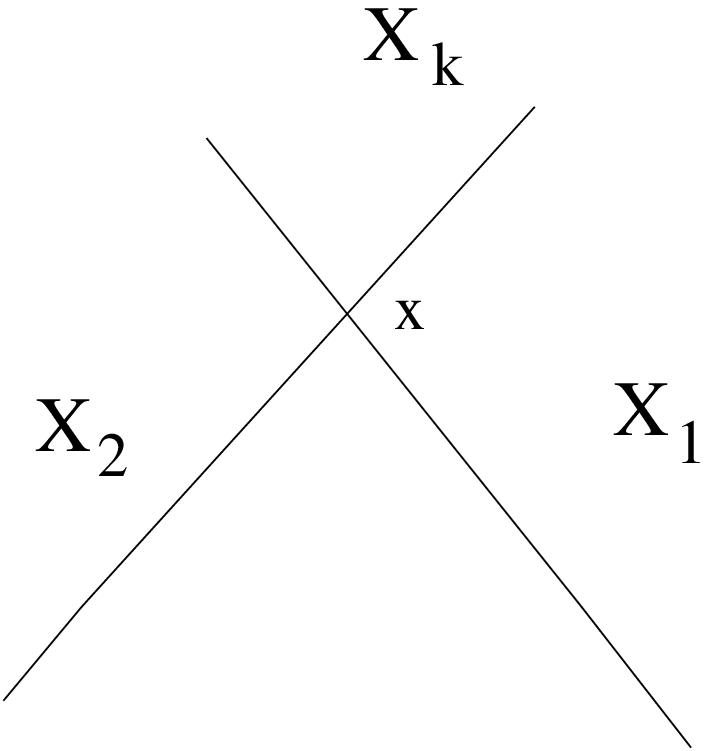}}

\pop {1}
\par
Consider a Galois cover $f:Y\to X$ of degree $p$ 
with $Y$ normal and such that
the Galois cover $f_K:Y_K\to X_K$ induced by $f$ above the generic fibre
is a $\mu_p$-torsor (hence in particular is \'etale). We assume that the 
special fibre $Y_k:=Y\times _Rk$ of $Y$ is reduced. For $i\in\{1,2\}$,
let $f_i:V_i\to U_i:=X_i-\{x\}$ be the torsor above $U_i$ induced by the 
cover $f$, and let $m_i$ be the conductor 
of $f_i$ at the point $x$. Then only the following two cases can ossur:
\par
Case 1)\ Either $f_i$ is an \'etale torsor in which case either $f_i$ ramify
above $x$ or $f_i$ extends to an \'etale torsor above $X_i$. In both cases
$V_i$ is smooth.
\par
Case 2)\ Or $f_i$ is a torsor under $\mu_p$ in which case is associated 
to $f_i$ canonically a logarithmic differential form $\omega _i$, and the 
singularities of $V_i$ lie above the zeros of $\omega_i$. In this case 
where $X_i$ is generic it was shown in  [Ra] (proposition 4) that 
such a differential form $\omega _i$ has $g_i-1$ double zeros if $p=2$ and 
$2g_i-2$ simple zeros if $p\neq 2$. In the later case $m_i=-2$ or $m_i=-1$
depending on whether $x$ is a zero of $\omega _i$ or not.
\par
In what follows we assume for simplicity that $p\neq 2$. Supppose we are in 
case 2, and let $y_j$ be a point of $V_i$ above a zero $z_j$ of $\omega_i$
which is different from the point $x$. Then $V_i$ is singular at $y_j$, the 
singularity at $y_j$ is unibranche, and
the geometric genus of $y_j$ equals $(p-1)/2$. Let $\Deg (z_j)$ 
be the non split
$\overline k$-degeneration data associated to $z_j$ via the semi-stable 
reduction of $Y$. Then  $\Deg (z_j)$ is necessarily of type 
$(0,(\mu_p, -2,0))$, and consists of one projective line $P_i$ with a marked 
$\overline k$-point $x_i=x$ and an \'etale torsor $T_i\to P_i-\{x_i\}$ with
conductor $2$ at the point $x_i$.
\pop {1}
\par
Above the double point $x$ three situations can occur:
\par
Case a)\ Both $f_1$ and $f_2$ are \'etale in which case $f$ is necessarily 
\'etale above $x$.

\pop {1}
\epsfysize=3cm
\centerline{\epsfbox{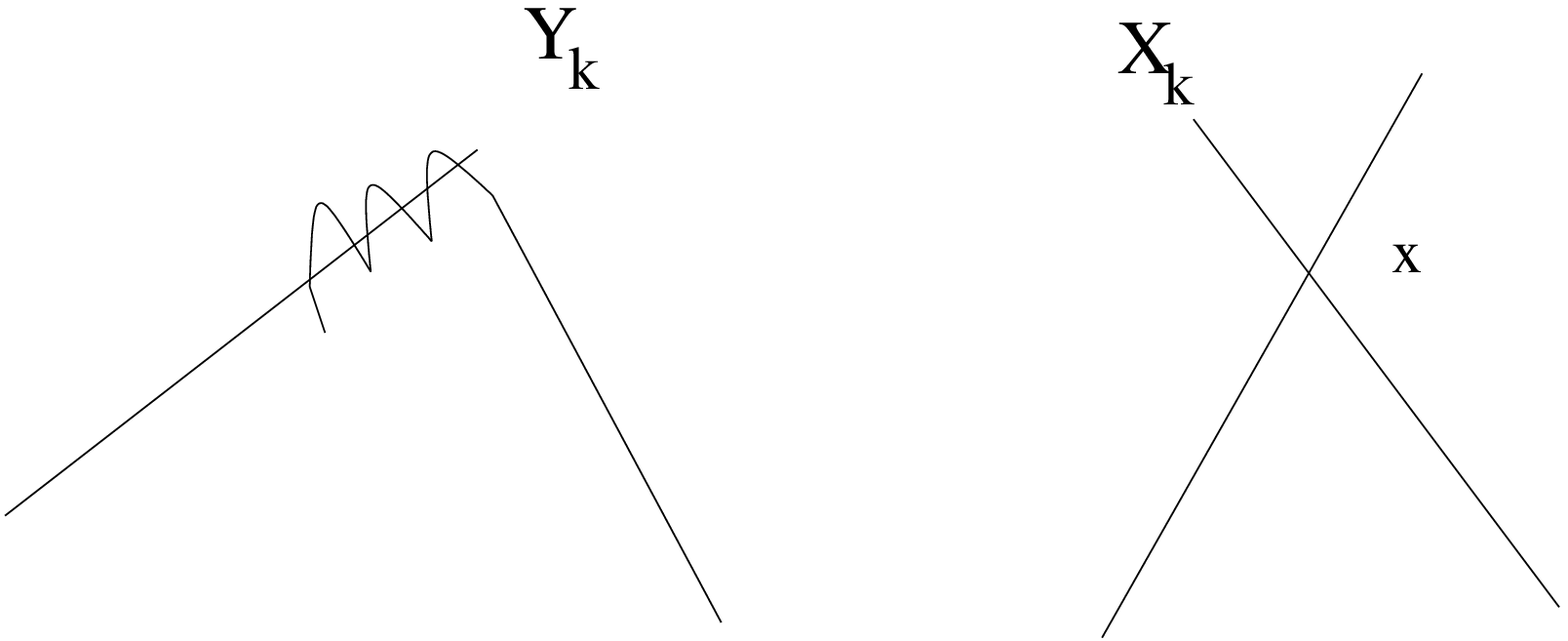}}

\pop {1}
\par
Apart from case a) $Y$ must be local above, we will 
denote by $y$ the unique point of $Y$ above $x$.

\par
Case b)\ Say $f_1$ is \'etale and $f_2$ is radicial, and assume for 
simplicity that $Y$ has two branches at the point $y$. In this case the 
geometric genus of $y$ equals $(m_1+2)(p-1)/2$ if $\omega _2$ has a zero at
$x$, or $(m_1+1)(p-1)/2$ otherwise.

\par
Case c)\ Both $f_1$ and $f_2$ are radicial in which case the geometric 
genus of $y$ equals $2(p-1)$, $3(p-1)/2$ or $p-1$ depending on whether
both $f_i$ have a zero at the point $x$, or just one of them do, or none of 
them.
\par
Let $\Deg (x)$ be the double non split $\overline k$-degeneration data
associated to $x$ via the semi-stable reduction of $Y$. assume for example 
that we are in case
c) and that the genus of $y$ equals $p-1$. Then the following 
two cases can occur:
\par
First case:\ either $\Deg (x)$ consists of one projective 
line $P$ with two marked 
$\overline k$-points $x_1=x$ and $x_2=x$, and an \'etale torsor 
$T\to P-\{x_1,x_2\}$ with conductor $2$ at each of $x_1$ and $x_2$. In this 
case the special fibre of the semi-stable reduction of $Y$ is a tree like
which consist of the components $V_1$ and $V_2$ which are linked by the 
projective completion of $V$ (which has genus $p-1$) at the points above $x$,
and at each of the points of $V_i$ above a zero of $\omega_i$ is linked a 
curve of genus $(p-1)/2$. 
\par
Second case:\ $\Deg (x)$ consists of two projective line $P_1$ 
and $P_2$ with a marked $\overline k$-points $x_1=x$ and $x_2=x$ on each one, 
and an \'etale torsor $T_i\to P_i-\{x_i\}$ with conductor $2$ at 
$x_i$. In this 
case the special fibre of the semi-stable reduction of $Y$ contains 
$(p-1)$ cycles and it consists of the components $V_1$ and $V_2$ which are 
linked by the projective completion of $T_1$ and $T_2$ (each has genus 
$(p-1)/2$) at the points above $x$, moreover both $T_1$ and $T_2$ meet at 
$p$ double points, and at each of the points of $V_i$ above a zero of 
$\omega_i$ is linked a curve of genus $(p-1)/2$ as in the preceeding case.

\pop {1}
\epsfysize=4cm
\centerline{\epsfbox{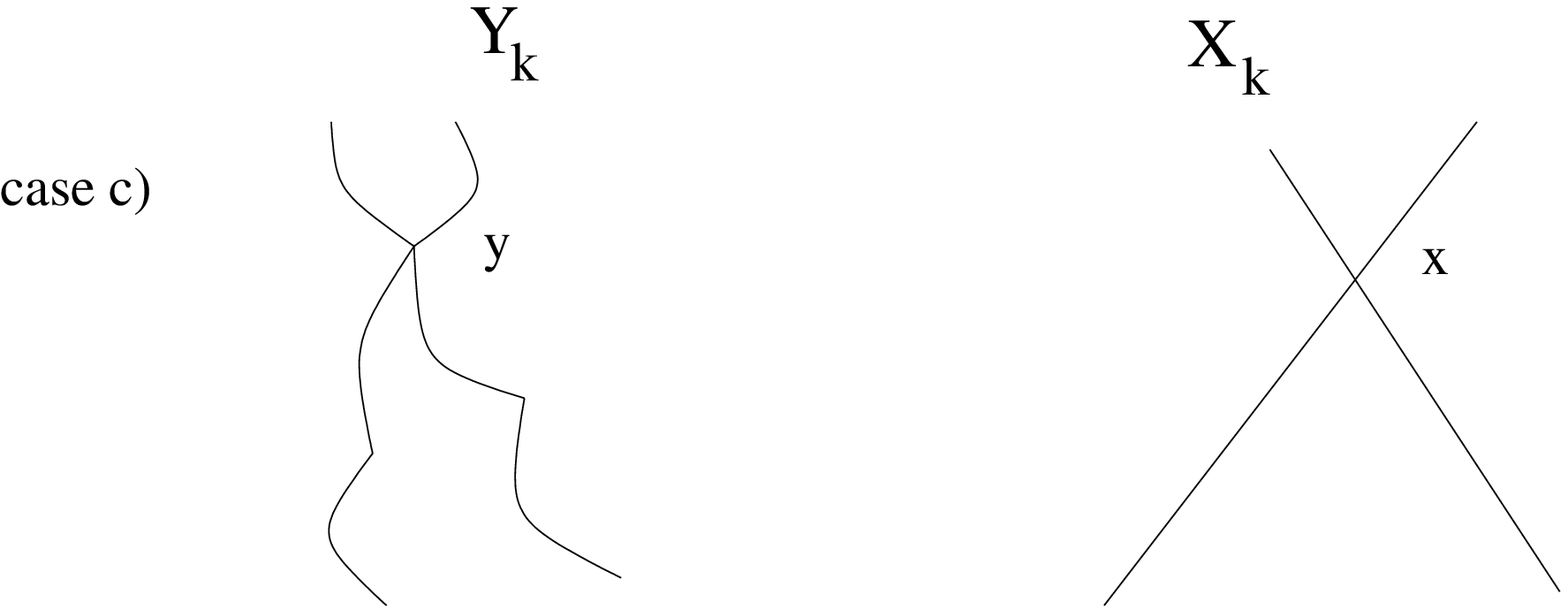}}

\pop {1}
\epsfysize=4cm
\centerline{\epsfbox{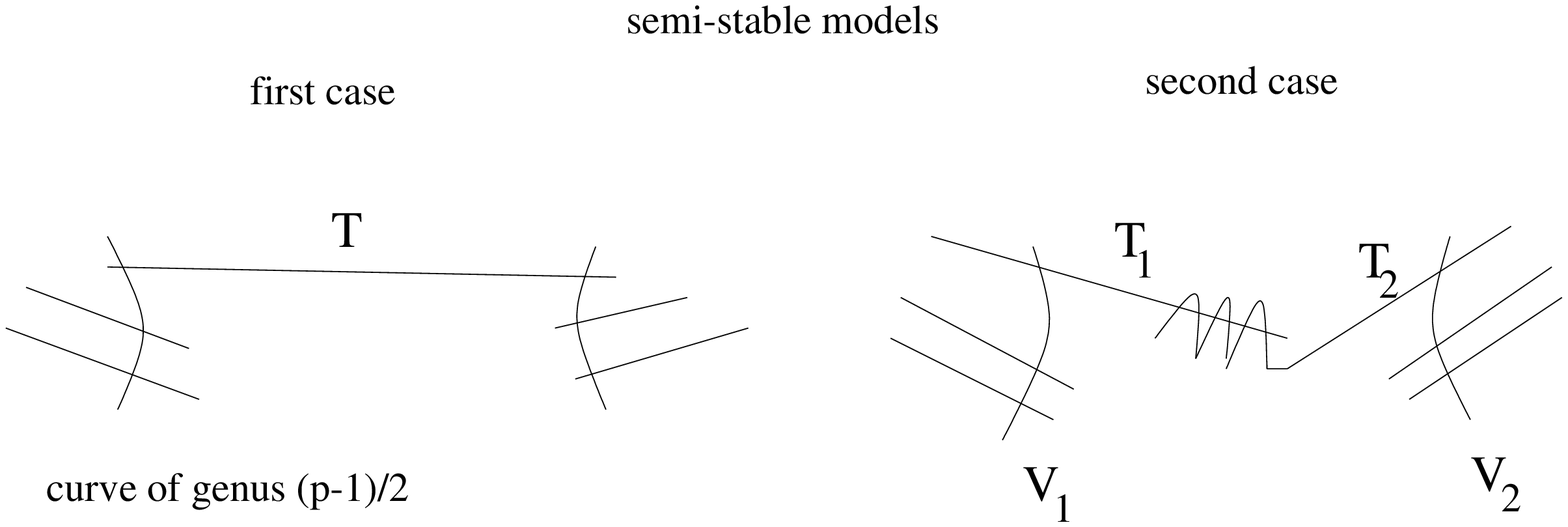}}
\pop {1}
\par
\noindent
{\bf \gr 2.3.}\rm \ There is a natural notion of isomorphism of 
degeneration data of rank $p$ 
associated to a proper and semi-stable $k$-curve $X_k$, over a given algebraic 
extension of $k$. We will denote by $\Deg_p(X_k)$ the set of 
isomorphism classes of $\overline k$-degeneration datas of rank $p$ 
associated to $X_k$. The set $\Deg_p(X_k)$ is equipped in a natural way with 
a $G_k$-action, via the natural action of $G_k$ on the above datas which was 
explained in [Sa] and I, and hence is a $G_k$-set. The following result is 
analogous to the result of realisation for degeneration datas obtained in 
paragraph I in the local situation.

\pop {.5}
\par
\noindent
{\bf \gr 2.3.1. Theorem.}\rm \ {\sl Let $X_k$ be a proper and semi-stable 
$k$-curve. Let $R$ be a complete discrete valuation ring of inequal 
characteristics with residue field
$k$, and fractions field $K$. Suppose given a 
$\overline k$-degeneration data $\deg (X_k)$ of rank $p$ 
associated to $X_k$ (in other words suppose given an element of $\Deg_p(X_k)$).
Then there exists, after eventually a finite extension of $R$, 
a proper and semi-stable $R$-curve $\Tilde X$ with smooth generic fibre 
and a special fibre isomorphic to $X_k$, and a Galois cover
$f:\Tilde Y\to \Tilde X$, such that the degeneration data 
associated to $f$ as in 2.2
is isomorphic to the given data $\deg (X_k)$. The above cover $f$ is not 
unique in general, moreover it can be constructed such that the cover $f_K:
\Tilde Y_K\to \Tilde X_K$ induced by $f$ on the 
generic fibres is ramified above 
$r=\sum _{t\in J}r_t+\sum _{i\in I}\sum _{j\in S_j}r_{i,j}$ 
geometric points of $X_K$ (hier $r_t$ and $r_{i,j}$ are the one given 
in the definition of the degeneration data in 2.2.1).}

\pop {.5}
\par
The  next result expreses the compatibility of degeneration with 
the Galois action.

\pop {.5}
\par
\noindent
{\bf \gr 2.3.2. Theorem.}\rm \ {\sl Let $R$ be a complete discrete valuation 
ring of inequal characteristics with residue field $k$, 
and fractions field $K$, and let $X$ be a proper 
and semi-stable $R$-curve with smooth and geometrically connected 
generic fibre $X_K$ and with 
special fibre $X_k$. Let $\overline K$ be an algebraic closure of $K$, and
let $L$ be the function field of the geometric generic fibre 
$X_{\overline K}:=X\times _R \overline K$ of $X$. Then the canonical 
specialisation map $\Sp: H^1_{\et}(\Spec L,\mu_p) \to \Deg _p(X_k)$ is 
$G_K$-equivariant, where hier we consider the natural action of $G_K$ on
$\Deg _p(X_k)$ via its canonical quotient $G_k$.}

\pop {.5}
\par
Note that the above specialisation map can not be surjective in this case
since the set $\Deg_p(X_k)$ with fixed ramification datas on the 
``generic fibre'' is infinite in general.
The proof of the above two results is very similar to the one obtained in 
paragraph I in the local case and are left to the reader. They are 
based on formal patching techniques and the results in [Sa] and [Sa-1]. 
However in order to prove 2.3.1 one will also need the following 
result on lifting of admissible torsors with given degeneration 
datas at the ``critical points''. More precisely let $R$ be as in 1.0. 
Let $X$ be a formal smooth $R$-curve, 
whose special fibre $X_k$ is irreducible.
Let $f_k:Y_k\to X_k$ be an admissible torsor, under a finite and 
flat $k$-group scheme
$G_k$ of rank $p$, which is radicial. Let $\omega$ be the associated 
differential form, and let $Z:=\{x_j\}_{j\in J}$ be the set 
of zeros of $\omega$ which are contained in $X_k$, which we call 
the {\bf critical} points of $f_k$, in particular the 
singularities of $Y_k$ lie above the points of $Z$. 
Assume that  $f:Y\to X$ is a torsor under a finite and 
flat $R$-group scheme $G$ which lifts $f_k$. To each point 
$x_j\in Z$ is associated via $f$, a simple degeneration data 
$\Deg (x_j)$ of type $(0,(G_k,m_j,h_j))$, where $m_j$ is the conductor of 
$f_k$ at $x_j$, and $h_j$ its residue at this point. The following result 
shows that we can choose such a lifting $f$ of $f_k$ which gives rise to a 
given set of degeneration datas at the critical points.

\pop {.5}
\par
\noindent
{\bf \gr 2.3.3. Proposition / Definition.}\rm \ {\sl Let $f_k:Y_k\to X_k$ be 
an admissible torsor under a finite and flat $k$-group scheme $G_k$ of 
rank $p$, which is radicial, and let $Z:=\{x_j\}_{j\in J}$ be the set 
of critical point of $f_k$ which we assume to be $k$-rational. For each point 
$x_j\in Z$ consider the pair 
$(m_j,h_j)$, where $m_j$ is the conductor of $f_k$ at $x_j$, and $h_j$ 
its residu at this point. Suppose given for each critical point $x_j$,
a simple degeneration data $\Deg (x_j)$ of type $(0,(G_k,m_j,h_j))$. 
Then there 
exists, after eventually a finite extension of $R$, a smooth formal $R$-curve 
$\Tilde X$ with a special fibre $\Tilde X_k:=\Tilde X\times _Rk$
which is isomorphic to $X_k$, and
a torsor $\Tilde f:\Tilde Y\to \Tilde X$ under a finite and flat $R$-group 
scheme $G$ of rank $p$, such that the restriction of $\Tilde f$ to $U_k:=X_k-Z$
is isomorphic to the restriction of $f_k$ to $U_k$, and such that the 
simple degeneration data assosiated to each point $x_j\in Z$,
via $\Tilde f$, is isomorphic to $\Deg (x_j)$. We call such a torsor 
$\Tilde f$ a {\bf lifting} of $f_k$,
with the {\bf given degenaration datas} $\{\Deg (x_j)\}_{j\in Z}$ at the 
critical points $\{x_j\}_{j\in Z}$.

\pop {.5}
\par
\noindent
{\bf \gr Proof.}\ \rm The torsor $f_k:Y_k\to X_k$ is admissible hence 
by definition can be lifted, after eventually a ramified extension of $R$, 
to a torsor
$f:Y'\to X$ under a finite and flat $R$-group scheme of rank $p$, which is 
either $\mu_p$ or $\Cal H_n$ for $0<n<v_K(\lambda)$ where $v_K$ is the 
normalised valuation of the fractions field of $R$ (such a lifting is 
not unique). Also for each point $x_j\in Z$ one can find a Galois cover
$f_j:\Cal Y_j\to \Cal X_j$ of degree $p$, where $\Cal X_j$ is the formal germ
of $X$ at $x_j$, and such that the degeneration data associated to $x_j$
via $f_j$ is isomorphic to $\Deg (x_j)$ (cf. 1.2.3). Let 
$U_k:=X_k-Z$, and let $g_k$ 
be the restriction of $f_k$ to $U_k$. Then $f$ induces a lifting $g:V\to U$ 
of $g_k$, which is a torsor under a finite and flat $R$-group scheme of 
rank $p$, and where $U$ is obtained from $X$ by deleting the formal fibres 
at the points $\{x_j\}_{j\in Z}$. Now a patching of the covers $g$ and 
$f_j$ along the formal fibres at the points $\{x_j\}_j$ will produce 
the desired torsor $\Tilde f:\Tilde Y\to \Tilde X$.

\pop {.5}
\par
\noindent
{\bf \gr 2.3.4. Example.}\rm \ Hier we use the same notations as in the 
example 2.2.3. First we will explain the Galois action on a $\mu_p$-torsor
$f_{\overline K}:Y_{\overline K}\to X_{\overline K}$ above the geometric 
generic fibre $X_{\overline K}$ of $X$. We assume that the above torsor 
induces in reduction the situation in case c) where $g_y=p-1$ and where the 
semi-stable reduction of $Y_{\overline K}$ is a tree like. Let 
$\sigma \in G_K$, and let $\overline \sigma$ be the image of $\sigma$ 
in $G_k$. Consider the $\mu_p$-torsor 
${f_{\overline K}}^{\sigma}:{Y_{\overline K}}^{\sigma}\to X_{\overline K}$
which is the conjugate of $f_{\overline K}$ via $\sigma$. Then using the 
result in 2.3.2 one can describe the semi-stable reduction of 
${Y_{\overline K}}^{\sigma}$. More precisely the graph of the semi-stable 
reduction of 
${Y_{\overline K}}^{\sigma}$ is isomorphic to the one of $Y_{\overline K}$, 
in particular it is also a tree like, and it consists of the components
$V_i^{\overline \sigma}$, for $i=1,2$, where ${V_i}^{\overline \sigma}$ is the 
conjugate of $V_i$ via ${\overline \sigma}$, which are linked by an irreducible
curve $T^{\overline \sigma}$ of genus $p-1$ which is the conjugate of $T$ by
${\overline \sigma}$. Moreover at each point of ${V_i}^{\overline \sigma}$
above a zero of $\omega_i$ is linked a curve of genus $(p-1)/2$ which is 
the conjugate by ${\overline \sigma}$ of the curve that is linked 
to the corresponding point (via the action of ${\overline \sigma}$) of 
$V_i$ in the graph of semi-stable reduction of $Y_{\overline K}$.

\par
Also one can establish the converse of the situation in 2.2.3, i.e. 
reconstruct a $\mu_p$-torsor $f_K$ as above from the degeneration datas, 
using 2.3.1.
More precisely, let $X_k$ be the semi-stable curve considered in 2.2.3 and
consider the non split degeneration data $\Deg (X_k)$ of rank $p$ which 
consists in the following: 
\par
For $i=1,2$ is given a $\mu_p$-torsor
$Y_i\to X_i$ associated to the logarithmic differential form $\omega_i$, 
and assume that the zeros $\{z_{j,i}\}_j$ of $\omega_i$ lie outside the point 
$x$. For each $i,j$ suppose given a non split simple degeneration data 
$\Deg (z_{j,i})$ of type $(0,(\mu_p,-2,0))$ and which consists of one 
projective line $P_{j,i}$ with one $\overline k$-marked point $z_{j,i}$ 
and an \'etale
torsor above $P_{j,i}-\{z_{j,i}\}$ with conductor $2$ at the point $z_{j,i}$.
At the double point $x$ is given a non split double degeneration data
of type $(0,(\mu_p,-2,0),(\mu_p,-2,0))$ and which consists of one projective 
line $P$ and two $\overline k$-marked points $x_1=x$ and $x_2=x$ and an 
\'etale torsor above $P-\{x_1,x_2\}$ with conductor $2$ at both 
points $x_1$ and $x_2$. Then by 2.3.1 one can find 
a proper and semi-stable curve $\Tilde X$ with smooth generic fibre
and a special fibre $\Tilde X_k:=\Tilde X\times R_k$ isomorphic to $X_k$,
and (eventually after a finite extension of $R$) a Galois cover 
$\Tilde f:\Tilde Y\to \Tilde X$ of degree $p$, such that the 
degeneration data of rank $p$ associated to $X_k$ via the semi-stable
reduction of $\Tilde Y$ is isomorphic to the above given data
$\Deg (X_k)$.

\pop {2}
\par
\noindent
{\bf \gr References}
\rm

\pop {.5}
\par
\noindent
[De-Mu] P. Deligne and D. Munford, {\sl The irreducibility of the space of 
curves with given genus}, Publ. Math. IHES. 36, 75-109, (1969). 

\pop {.5}
\par
\noindent
[Gr] A. Grothendieck, S\'eminaire de g\'eom\'etrie alg\'ebrique SGA-1, 
Lecture Notes 224, Springer Verlag, (1971).

\pop {.5}
\par
\noindent
[He] Y. Henrio, {\sl Arbres de Hurwitz et automorphismes d'ordre $p$ des 
disques et couronnes $p$-adiques formels}. Th\`ese de doctorat, 
Bordeaux France (1999).

\pop {.5}
\par
\noindent
[Ra] M. Raynaud, {\sl $p$-Groupes et reduction semi-stable des courbes},
The Grothendieck Festschrift, vol. 3, 179-197, Birkh\"auser, (1990).

\pop {.5}
\par
\noindent
[Sa] M. Sa\"\i di, {\sl Torosrs under finite and flat group schemes 
of rank $p$ with Galois action}, preprint.

\pop {.5}
\par
\noindent
[Sa-1] M. Sa\"\i di, {\sl Wild ramification and vanishing cycles formula},
preprint.

\pop {2}
Mohamed Sa\"\i di

\pop {1}
\par
Departement of Mathematics
\par
University of Durham
\par
Science Laboratories
\par
South Road
\par
Durham, DH1 3LE, UK
\par
saidi\@durham.ac.uk

\enddocument
\end